\documentclass[12pt,a4paper]{amsproc}

\usepackage[T1]{fontenc}
\usepackage[utf8]{inputenc}
\usepackage{lmodern}
\usepackage{amsmath,amssymb,amsthm,mathtools, mathrsfs}
\usepackage{enumitem, todonotes}
\usepackage[
  colorlinks=true,
  linkcolor=blue,
  citecolor=blue,
  urlcolor=blue,
  pagebackref=true
]{hyperref}

\usepackage[
  a4paper,
  left=29mm,
  right=29mm,
  top=29mm,
  bottom=31mm
]{geometry}

\usepackage{tikz}

\title[Metric structure in Hecke-induced modules]
{Metric structure in Hecke-induced modules\\ for  groups acting on buildings}

\author[P.-E. Caprace]{Pierre-Emmanuel Caprace}
\address{P.-E.C., UCLouvain, Chemin du Cyclotron 2, 1348 Louvain-la-Neuve, Belgium}
\email{pierre-emmanuel.caprace@uclouvain.be}

\author[A. Lou\'e]{Alex Lou\'e}
\address{A.L., UCLouvain, Chemin du Cyclotron 2, 1348 Louvain-la-Neuve, Belgium}
\email{alex.loue@uclouvain.be}

\author[A. Thom]{Andreas Thom}
\address{A.T., TU Dresden, 01062 Dresden, Germany}
\email{andreas.thom@tu-dresden.de}
\date{November 23, 2026}

\thanks{PEC is supported in part by the FWO and the F.R.S.-FNRS under the EOS programme (project ID 40007542).}

\newtheorem{theorem}{Theorem}[section]
\newtheorem{thmintro}{Theorem}

\newtheorem{corintro}[thmintro]{Corollary}
\newtheorem{proposition}[theorem]{Proposition}
\newtheorem{lemma}[theorem]{Lemma}
\newtheorem{corollary}[theorem]{Corollary}

\theoremstyle{definition}

\newtheorem{remark}[theorem]{Remark}
\newtheorem{example}[theorem]{Example}

\DeclareMathOperator{\Aut}{Aut}
\DeclareMathOperator{\End}{End}
\DeclareMathOperator{\Ind}{ind}
\DeclareMathOperator{\Stab}{Stab}
\DeclareMathOperator{\Rad}{Rad}

\DeclareMathOperator{\supp}{supp}

\DeclareMathOperator{\Ch}{Ch}

\DeclareMathOperator{\proj}{proj}
\DeclareMathOperator{\Res}{Res}
\DeclareMathOperator{\dist}{dist}

\newcommand{\CC}{\mathbf C}
\newcommand{\KK}{\mathbf K}
\newcommand{\RR}{\mathbf R}

\newcommand{\one}{\mathbf 1}

\newcommand{\qq}{\mathbf q}
\newcommand{\rg}{\mathrm g}

\begin{document}

\begin{abstract}
Let $X$ be a thick locally finite building of type $(W, S)$ and $G$ be a locally compact group acting properly and chamber-transitively on $X$. We give a condition on a Hermitian module $(V_0, \varphi_0)$ of the Iwahori--Hecke algebra $\CC_\qq[W]$ ensuring that the natural $G$-invariant form defined on the induced $G$-module has the same negative index as $\varphi_0$, where $\qq$ is the thickness vector of $X$. It implies that induction of Hecke-modules preserves unitarity. In the affine case, this recovers results by Matsumoto, Barbasch--Moy and Ciubotaru. In general, assuming that $G$ is Weyl-transitive, we obtain new irreducible $G$-representations  on Hilbert spaces as well as on infinite-dimensional symmetric spaces of finite rank. Using Hecke-modules of    Lorentzian signature, we also derive new obstructions to Kazhdan's property (T), showing that the known sufficient conditions in terms of the thickness are mostly necessary. 
\end{abstract}

\maketitle

\section{Introduction}
\label{sec:introduction}

Let $G$ be a totally disconnected locally compact group and $B \leq G$ be a compact open subgroup. The   algebra $C_c(B\backslash G/B)$ of bi-$B$-invariant functions on $G$ is a sub-$*$-algebra of the convolution algebra $C_c(G)$ of continuous, compactly supported, complex valued functions on $G$. It can be identified with the commutant of the linear representation of $G$ on the space $C_c(G/B)$, see \cite[Ch.~IV, Exer.~22]{Bour}. 

The  convolution algebra $C_c(B\backslash G/B)$  occurs    naturally in the representation theory of $G$. Indeed,   every unitary representation $(\mathcal H, \pi)$ of $G$ gives rise to a $*$-representation of $C_c(B\backslash G/B)$ on  the space $\mathcal H^B$ of $B$-invariant vectors. Moreover, if $\pi$ is irreducible and if $\mathcal H^B$ is non-zero, then $\mathcal H^B$ is an irreducible $C_c(B\backslash G/B)$-module, see \cite[Ch.~IV, Theorem 2]{Lang}. It is a natural question to determine whether every (irreducible) unitary $C_c(B\backslash G/B)$-module arises in this way. 

In general, the answer is negative:  counterexamples may be found with $G$  a simple algebraic group over a local field and $B$ is a \textit{maximal compact subgroup}, see \cite{LarsenPalma} and references therein. However, for a simple algebraic group $G$ over a local field, it turns out that the answer is positive if $B$ is an \textit{Iwahori subgroup}. This was established in the rank~one case by Matsumoto~\cite[\S 5.6.6, p.~215]{Matsumoto}. The question whether it also holds in higher rank was raised by Borel~\cite[\S 5.11]{Borel}, and remained a long standing problem until an affirmative answer was obtained by Barbasch--Moy and Ciubotaru~\cite{BM89, BM93, Ciubotaru}. 
 
It is natural to extend the scope of Borel's question to the case where $G$ acts properly, continuously and Weyl-transitively on a thick locally finite building $X$ of arbitrary type $(W, S)$, and $B$ is the stabilizer of a chamber. Beyond the affine case, explicit families of examples  are afforded by the maximal Kac--Moody groups over finite fields, in which case $(W, S)$ can be an arbitrary crystallographic Coxeter system  (see \cite[Ch.~8]{Marquis_book}). Moreover, in the right-angled case, additional examples that are not of Kac--Moody type arise when considering  the full automorphism group of semi-regular right-angled buildings (see \cite{CapRAB}), as well as variants with prescribed local actions (see \cite{DMSS}).  

In this case, the convolution algebra  $C_c(B\backslash G/B)$ is isomorphic to the Hecke algebra $\CC_{\mathbf q}[W]$   associated with the thickness vector $\mathbf q = (q_s)_{s \in S}$  of $X$ (see Section~\ref{sec:Hecke} below). Given a unitary $\CC_{\mathbf q}[W]$-module $V_0$, the \textbf{Hecke-induced module} 
$$C_c(G) \otimes_{C_c(B\backslash G/B)} V_0$$ 
carries a natural $G$-invariant form, called the \textbf{induced form}, see Section~\ref{sec:Hecke} below. Borel's question amounts to asking whether that induced form is positive semi-definite. The following result provides a positive answer for all types. 

\begin{thmintro}\label{thmintro:unitary}
Let $X$ be a thick locally finite building of type $(W, S)$ and $G$ be a locally compact group acting properly, continuously and Weyl-transitively on $X$. Let  $\mathbf q = (q_s)_{s\in S}$ be the thickness vector of $X$. For every unitary $\CC_{\mathbf q}[W]$-module,  the Hecke-induced form is positive semi-definite. 
\end{thmintro}

A unitary representation  $(\mathcal H, \pi)$ of $G$ is called \textbf{$B$-spherical} if the space $\mathcal H^B$ of $B$-invariant vectors is non-zero. 

\begin{corintro}\label{cor:equivalence}
Retain the hypotheses of Theorem~\ref{thmintro:unitary}. The map $\mathcal H \mapsto \mathcal H^B$ that associates to a  unitary representation $(\mathcal H, \pi)$ of $G$  the corresponding $C_c(B\backslash G/B)$-module $ \mathcal H^B$ induces a bijective correspondence between   equivalence classes of irreducible $B$-spherical representations of $G$ and the irreducible unitary $C_c(B\backslash G/B)$-modules.
\end{corintro}

Our proof is direct and self-contained; in particular we obtain a new proof in the affine case, significantly different from the aforementioned works \cite{BM89, BM93, Ciubotaru, Matsumoto}. It is uniform, in the sense that it does not distinguish between the principal series and the complementary series. An important ingredient is an exhaustion of the chambers of a building that satisfies some weak convexity property   (see Proposition~\ref{prop:exhaustion}), and whose construction exploits the interplay between the CAT($0$)  and the combinatorial metrics.

\medskip

The quadratic relations satisfied by the generators $(e_s)_{s \in S}$ of $\CC_{\mathbf q}[W]$, recalled in \S\ref{sec:hecke-preliminaries} below, ensure that the Hecke algebra	$\CC_{\mathbf q}[W]$ admits an enveloping $C^*$-algebra, which is denoted by $C^*_{\mathbf q}(W)$, see \cite[\S 2.4.2]{Matsumoto}. Theorem~\ref{thmintro:unitary} has the following consequence, where  $p_B = \frac 1 {m(B)} \one_B \in C^*(G)$ denotes the normalized characteristic function of $B$.

\begin{corintro}\label{corintro:corner}
Retain the hypotheses of Theorem~\ref{thmintro:unitary}, and normalize the Haar measure $m$ on $G$ so that $m(B)=1$. The canonical surjective map
$$C^*_{\mathbf q}(W) \to p_B C^*(G) p_B,$$ 
defined by the assignments $e_w \mapsto \one_{BwB}$, is injective. In particular $C^*_{\mathbf q}(W) $ and the corner algebra $p_B C^*(G) p_B$ are isomorphic.  
\end{corintro}
	
In case $(W, S)$ is of right-angled type, the Hecke algebra 	$\CC_{\mathbf q}[W]$ is isomorphic to the group algebra $\CC[W]$ (see \cite[Proposition~4.2]{CKL}), so that the enveloping $C^*$-algebra $C^*_{\mathbf q}(W)$ is isomorphic to the maximal $C^*$-algebra $C^*(W)$. No such description of $C^*_{\mathbf q}(W)$ is available in general. 
	
 Every one-dimensional representation of $\CC_{\mathbf q}[W]$ is unitary. Among them, the \textbf{alternating representation} $\CC_{\mathbf q}[W]\to \CC$ is the representation that sends the generator $e_s$ to $-1$ for all $s \in S$. In \cite[\S 1]{MatsumotoCRAS}, Matsumoto defines the \textbf{special representation} (also called the \textbf{Steinberg representation}) of $G$ as an irreducible unitary representation $(\mathcal H, \sigma)$  of $G$ such that $\mathcal H^B$ is one-dimensional, and the $\CC_{\mathbf q}[W]$-representation on $\mathcal H^B$ induced by $\sigma$ is the alternating representation. Matsumoto notes that,  if the special representation exists, then it is unique up to equivalence. Moreover,   the special representation exists provided $\sum_{w \in W}  q_w^{-1} < \infty$, where $q_w = m(BwB)$. In that case, $\sigma$ is square-integrable. (Those  assertions are made without proof in \cite{MatsumotoCRAS}; they can be established using the idempotents of the Hecke-von Neumann algebras constructed in \cite[\S 19.2]{Davis}. Another construction of the special representation in the square-integrable case appears in \cite[\S 6]{DymJan}.) The convergence of the series  $\sum_{w \in W}q_w^{-1} $ is automatic if $X$ is a thick building of affine type, but it fails in general, see \cite{Bischof}, so that the existence of the special representation has remained undecided in the corresponding cases.   Theorem~\ref{thmintro:unitary} has the following direct consequence. 
 
 \begin{corintro}	[Existence of the special representation]\label{corintro:special}			
Every group satisfying   the hypotheses of Theorem~\ref{thmintro:unitary} has a special representation. 
 \end{corintro}
 
 There is another  natural finite-dimensional representation $V_{\mathrm{g}}$ of $\CC_{\mathbf q}[W]$, in dimension $|S|$, that we call the \textbf{geometric Hecke module}, see Section~\ref{sec:geom-Hecke-module} below. It can be viewed as a deformation of the  geometric representation of the Coxeter group $W$. In case $(W, S)$ is \textbf{$2$-spherical} (i.e. the Coxeter number $m_{s, t}$ is finite for all $s \neq t \in S$), it carries a  canonical Hermitian form $\varphi_{\mathrm{g}}$, which is  the natural deformation of the Tits form defined on the geometric representation of $W$. That form $\varphi_{\mathrm{g}}$ need not be positive definite. But when it is, the Hecke-induced unitary representation of $G$ afforded by Theorem~\ref{thmintro:unitary} is called the \textbf{geometric unitary representation} of $G$. Using Kassabov's criterion for Kazhdan's property (T) (see \cite{Ka11}), we establish the following. 
 
 \begin{corintro}[Spectral gap of the Laplacian]\label{corintro:Kazhdan}
Retain the hypotheses of Theorem~\ref{thmintro:unitary} and suppose that $(W, S)$ is $2$-spherical. 
If $\varphi_{\mathrm{g}}$ is positive definite, then $G$ has Kazhdan's property (T). Moreover, the spectral gap of the Laplacian defined by 
$$\Delta = \sum_{s \in S} (\one_B - (q_s+1)^{-1} \one_{P_s}) \in C^*(G), $$
where $(P_s)_{s \in S}$ denote the rank one parabolic subgroups, is achieved by the geometric unitary representation, which is irreducible provided $(W, S)$ is irreducible. 
 \end{corintro}

The condition that  $(W, S)$ be $2$-spherical is well known to be a necessary condition for $G$ to have property (T), see \cite[\S2]{DymJan}.  We refer to  Corollary~\ref{cor:Kazhdan} for a more general statement on property~(T) for groups acting transitively (but not necessarily Weyl-transitively). 

Corollary~\ref{corintro:Kazhdan} can be viewed as a ``thick analogue'' of Kassabov's result  \cite[\S 6.1]{Ka11} that the geometric representation of a finite Coxeter group achieves the spectral gap of the Laplacian for  that group with respect to the standard Coxeter generating set. 

\medskip
Our method for proving Theorem~\ref{thmintro:unitary} applies not only to unitary $\CC_{\mathbf q}[W]$-modules, but also to  Hermitian  $\CC_{\mathbf q}[W]$-modules whose associated form has a non-trivial signature. In fact, the proof does not require that the building $X$ admits a Weyl-transitive automorphism group, or even a chamber-transitive automorphism group. Indeed, our main result is the following. 

\begin{thmintro}\label{thmintro:signature}
Let $X$ be a thick, locally finite, semi-regular, building of type $(W, S)$, all of whose spherical residues have a transitive automorphism group.  Let $\mathbf q = (q_s)_{s\in S}$ be its thickness vector, and let $(V_0, \varphi_0)$ be a non-degenerate Hermitian $\KK_{\mathbf q}[W]$-module, where $\KK$ is the field of real or complex numbers. 
We suppose that the form $\varphi_0$ has negative index $p< \infty$.  If for each spherical subset $J  \subset S$, there is a   $p$-dimensional subspace $U_J$ of $V_0$ that is pointwise fixed by $T_s =  q_s^{-1} e_s$  for all $s \in J$, and on which the restriction of $\varphi_0$ is negative definite, then  the  induced form on the Hecke-induced module $\Ind(V_0)$ also has negative index~$p$. 
\end{thmintro}

 A building is called is \textbf{semi-regular} if  for each $s \in S$, all $s$-panels of $X$ have the same cardinality. This is the condition needed for the thickness vector to be well defined. Observe that if $\Aut(X)$ is chamber-transitive, then $X$ is semi-regular and all spherical residues have a transitive automorphism group.

In case $\varphi_0$ is positive definite, we have $p =0$, and the  assertion on spherical subsets  $J  \subset S$  is redundant. Hence Theorem~\ref{thmintro:unitary} follows as a direct consequence of Theorem~\ref{thmintro:signature}. 
 
The   assertion  on spherical subsets in Theorem~\ref{thmintro:signature} is automatically satisfied by the geometric Hecke module  (see Lemma~\ref{lem:geom-form-tech} below). Using that the existence of an unbounded representation preserving a form with signature $(1, \infty)$ is an obstruction to property (T), we   deduce the following result, that may be viewed as a partial converse to Corollary~\ref{corintro:Kazhdan}.

\begin{corintro}\label{corintro:non-Kazhdan}
Let $X$ be a thick locally finite building of $2$-spherical type $(W, S)$. Suppose that  $S$ contains a spherical or affine subset of size~$|S|-1$, and that the geometric form $\varphi_\rg$ is non-degenerate. Then  a locally compact group $G$ acting properly, continuously and chamber-transitively on $X$ has Kazhdan's property (T) if and only if $\varphi_{\mathrm{g}}$ is positive definite.
\end{corintro}

We refer to  Corollary~\ref{cor:Lorentzian} for a  more general criterion ensuring the failure of property~(T).

The first results establishing property (T) for groups acting strongly transitively on thick locally finite buildings of $2$-spherical type and sufficiently large thickness, beyond the affine case, were due to Dymara--Januszkiewicz \cite{DymJan}. Their results required a lower bound on the thickness that could subsequently be improved by invoking results of Ershov--Jaikin~\cite{EJ},  Kassabov~\cite{Ka11}  and Ershov--Rall~\cite{ErshovRall} (see also \cite{GriOpp} for a reformulation of the Kassabov threshold in that context). The status of property (T) for groups acting on $2$-spherical buildings of low thickness has long remained unclear. For buildings of affine type, a recent result of I.~Oppenheim~\cite{Oppenheim} ensures that property (T) holds unconditionally. In the non-affine case, to the best of our knowledge, the only previously known result providing an obstruction to property (T) in that context  so far was Bischof's negative result \cite[Lemma~7.4.6]{BischofThesis} on Kac--Moody hyperbolic buildings of type $(4, 4, 4)$ over $\mathbf F_2$.
\medskip 

For buildings of compact hyperbolic type whose geometric form is non-degenerate, the alternative in Corollary~\ref{corintro:non-Kazhdan} can be strengthened to Kazhdan versus Haagerup property.
Recall that a Coxeter system is of \textbf{compact hyperbolic type}
if it is non-spherical, non-affine and if for every $J \subsetneq S$, the parabolic subgroup $W_J = \langle J \rangle$ is finite.  Note that in this case, the geometric form $\varphi_\rg$ is either positive definite, degenerate positive semi-definite, or non-degenerate of Lorentzian signature (see Corollary~\ref{cor:lorentzian-dichotomy}). Using a case-by-case computation, one may check that  the form $\varphi_\rg$ is non-degenerate if the thickness vector is constant.

\begin{corintro}
\label{corintro:Haagerup}
Let $X$ be a thick locally finite building of compact hyperbolic type $(W,S)$, of rank at least three, and let $\qq=(q_s)_{s\in S}$ be its thickness vector. Let $G$ be a locally compact group acting properly, continuously and chamber-transitively by automorphisms on $X$.

If $\varphi_\rg$ has Lorentzian signature, then $G$ has the Haagerup property, and the hyperbolic action obtained from the geometric Hecke module is metrically proper.
\end{corintro}

This recovers the special case of Kac--Moody  buildings of type $(4, 4, 4)$  over the field of order~$2$, recorded in \cite[Proposition~7.1]{CapTho}. Corollary~\ref{corintro:Haagerup} will be deduced from a more general statement (see Proposition~\ref{prop:proper}) affording the Haagerup property under more flexible hypotheses, that also applies beyond the compact hyperbolic type, see \S\ref{sec:proper}. 
	
\medskip	
 In case $(W, S)$ is right-angled and $G$ is strongly transitive on $X$ (this is automatically the case if $G = \Aut(X)$ and $X$ is semi-regular, see \cite{CapRAB}), the geometric Hecke module is not canonically defined: it depends on a choice of parameters $\mu_{s, t} \in \CC$ for all pairs $s \neq t \in S$ for which the Coxeter number $m_{s, t}$ is infinite. Using Theorem~\ref{thm:action-sym-space}, we obtain a multi-parameter family of irreducible representations of those groups on infinite-dimensional symmetric spaces of finite rank. In the special case where $X$ is a semi-regular tree and $\mu_{s, t} >1$ is real, we recover the irreducible representations on the infinite-dimensional real hyperbolic space first constructed and studied by Burger--Iozzi--Monod~\cite{BIM}. In our context, irreducible representations on infinite-dimensional symmetric spaces of   finite  rank~$>1$ also occur naturally, not only in the right-angled case but also for some $2$-spherical types in low thickness, see \S\ref{sec:examples} below for   explicit examples. This contrasts with the classical case of simple Lie and algebraic groups of higher rank, for which various results ensure that such actions cannot exist, see the recent work of B.~Duchesne \cite{Duchesne_2023} and F.~Viola \cite{Viola}.

\subsection*{Acknowledgements}

We thank Bruno Duchesne, Timothée Marquis, Nicolas Monod, Pierre Py, Sven Raum and Alain Valette for their comments on this work, and Bruno Duchesne for drawing our attention to the work of F.~Viola~\cite{Viola}.

\subsection*{AI statement}
The idea that a control of the signature of the induced form on a Hecke-induced module would yield an answer to Borel's question for locally compact groups with BN-pairs of general type, as well as an obstruction to Kazhdan's property (T), was   humanly designed several years ago by the first two authors. This idea was shared with the third author upon the completion of \cite{CapTho}, when the question, whether the Haagerup property holds for groups acting on buildings whose thickness is below the Kassabov threshold, naturally occurred. An intermediate result on that specific question,  based on a combinatorial strategy implementing that idea,  was then obtained  by the last author in July 2026, relying on human--machine interactions with GPT 5.6 Sol.  Its target was the Haagerup property  for groups acting Weyl-transitively on buildings of compact hyperbolic type, constant thickness and rank~$3$, via the Hecke-induced geometric representation. This has led us to the idea of the base change in Proposition~\ref{prop:add-residue}. The paper was then developed collaboratively by the three authors, without   further use of any AI  tool, until a complete draft without \S\ref{sec:proper} had been written. In a final phase, an interaction with  ChatGPT 6 Astra helped in obtaining Corollary~\ref{corintro:Haagerup}, that we then generalized to Proposition~\ref{prop:proper}. The whole document was lastly tested with   ChatGPT 6 Astra for consistency, which led us to make a few minor corrections.  

\section{Buildings}

\subsection{Combinatorial structure and  CAT(0) realization}

For the general theory, we refer to the books \cite{AbramenkoBrown} and \cite{Davis}. We briefly recall some  basic definitions and facts that are needed for our purposes. 

Let $(W,S)$ be a Coxeter system. We assume throughout that $S$ is finite. We shall view a building $X$ of type $(W, S)$ as a set of chambers, denoted by $\Ch(X)$, endowed with a $W$-valued map $\delta \colon \Ch(X) \times \Ch(X) \to W$ called the \textbf{Weyl-distance}, as in \cite[Chapter~5]{AbramenkoBrown}. An \textbf{automorphism} of the building $X$ is a permutation of $\Ch(X)$ that preserves the Weyl-distance. (Note that automorphisms are  type-preserving by definition in this approach to buildings.) Given $c, d \in \Ch(X)$, the word length of $\delta(c, d)$ with respect to the Coxeter generating set $S$ is called the \textbf{numerical distance} (or \textbf{gallery distance}) between $c$ and $d$. It is denoted by $\dist(c, d)$. Given a residue $R$ in $X$, we denote by $\Ch(R)$ its set of chambers. We say that are residue $R_1$ is contained in a residue $R_2$ if $\Ch(R_1) \subseteq \Ch(R_2)$. We say that a residue $R$ is contained in an apartment $A$ if some chamber of $R$ belongs to $A$.

The Davis realization of $X$  is denoted by $|X|$. It is defined  as follows (see \cite[\S18]{Davis} for more details). Consider  the poset $\mathcal S(X)$ formed by the spherical residues in $X$. The \textbf{flag complex} of the poset $\mathcal S(X)$, denoted by $\mathrm{Flag}(\mathcal S(X))$, is the abstract simplicial complex whose simplices are chains of spherical residues. In particular, vertices are spherical residues of $X$. The space $|X|$ is then defined as the geometric realization of that simplicial complex $\mathrm{Flag}(\mathcal S(X))$. The space $|X|$ is endowed with the CW topology (see \cite[\S A.1]{Davis}). A crucial feature is that it also carries a piecewise Euclidean structure such that $|X|$, once endowed with  the associated length metric denoted by $d$,    is  a complete  CAT($0$) metric space, see \cite[Theorem~18.3.1]{Davis}.  

Given a chamber $c \in \Ch(X)$, we denote by $|c|$ the corresponding subcomplex of $|X|$. Its simplices are the chains of spherical residues that contain the chamber $c$. It is also a  closed convex subset of the CAT($0$) space $|X|$, which is moreover compact since $S$ is finite. 

Given a spherical residue $R$, we set 
$$|R| = \bigcap_{c \in \Ch(R)} |c|$$
and
$$|\Ch(R)| = \bigcup_{c \in \Ch(R)} |c|.$$

Clearly $|R|$ is closed and convex. 

\begin{lemma}\label{lem:Ch(R)}
For any spherical residue $R$, the set $|\Ch(R)|$ is closed and convex. 
\end{lemma}
\begin{proof}
Given two chambers $c, d \in \Ch(R)$, consider an apartment $A$ containing them both. The intersection $\Ch(R) \cap A$ coincides with the intersection of all half-apartments that contain $c$ and $d$. Since the half-apartments are convex in the CAT($0$) realization, it follows that $|\Ch(R)| \cap |A|$ is convex. Hence the geodesic segment joining two points of $|c| \cup |d|$ is contained in $|\Ch(R)|$. Since this holds for any pair of chambers of $R$, we infer that $|\Ch(R)|$ is indeed convex. 

By construction $|\Ch(R)|$ is a piecewise Euclidean complex with finitely many isometry types of cells. Hence it is Cauchy complete (see \cite[pp.~105--111]{BH99} or \cite[Proposition~I.3.4]{Davis}), and therefore   a closed subset of $|X|$. 
\end{proof}

Each $x \in |X|$ has a \textbf{support}, which is the unique simplex $\sigma \in \mathrm{Flag}(\mathcal S(X))$ that contains $x$ in its relative interior. Recall that, by definition, $\sigma$ is a chain of spherical residues in $X$ ; we denote by $\Res(x)$ the smallest spherical residue that belongs to $\sigma$.
We set $\Ch(x) = \Ch(\Res(x))$ and $|\Ch(x)| = |\Ch(\Res(x))|$.

\begin{lemma}\label{lem:Ch(x)}
For all $x \in |X|$, the set $|\Ch(x)|$ is a neighborhood of $x$ in $|X|$.
\end{lemma}
\begin{proof}
Let $\sigma\in \mathrm{Flag}(\mathcal S(X))$ be the support of $x$. Let also $\tau \in \mathrm{Flag}(\mathcal S(X))$ be a simplex  containing $\sigma$. Then the smallest spherical residue $T\in \tau$ is contained in the smallest spherical residue of $\sigma$, which is $\Res(x)$ by definition. Thus any chamber $c \in \Ch(T)$ belongs to $\Ch(\Res(x))$. Since $|c|$ is a subomplex of $|X|$ that contains the simplex $\tau$, we infer that $|\Ch(x)|$ contains $\tau$. In particular $|\Ch(x)|$ contains the union $\mathrm{St(\sigma)}$ of all simplices $\tau$ with $\tau \supseteq \sigma$. That set  is a neighborhood of $x$.
\end{proof}

\subsection{Projections}

A fundamental feature of buildings is the existence of combinatorial projections between residues. We present basic facts and refer to  \cite[\S 3.19]{Tits74} and \cite[\S2]{CapRAB} for more information.

Given a chamber $c \in \Ch(X)$ and a residue $R$ in $X$, the \textbf{projection} of $c$ on $R$ is the unique chamber of $\Ch(R)$ that is closest to $c$ in the numerical distance. 
It is denoted by $\proj_R(c)$. An important property of the projection map is that it does not increase the numerical distance between chambers: for all $c, c' \in \Ch(X)$, the numerical distance from $\proj_R(c)$ to $\proj_R(c')$ is bounded above by the numerical distance from $c$ to $c'$. 

If $R$ and $R'$ are two residues, then the set 
$$\{\proj_R(c) \; | \; c \in \Ch(R')\}$$
is the chamber-set of a residue contained in $R$. That residue is denoted by $\proj_R(R')$. 

The following result provides relevant information on $\proj_R(R')$ in the case of thin buildings. This is particularly useful, since $\proj_R(R')$ is completely determined by the corresponding projection taken in any apartment meeting $R$ and $R'$ (see \cite[\S 3.19]{Tits74}). Given a residue $R$ of a thin building $X$ of type $(W, S)$, we denote by $\mathcal H(R) $ be the collection of all hyperplanes of  $X$ that cross   $\Ch(R)$. Viewing $W$ as the automorphism group of $X$, we denote by $W_R$ the stabilizer of $R$ in $W$. It is a parabolic subgroup of $W$. 

\begin{lemma}\label{lem:proj_R(R')}
Let $R, R'$ be two residues in a thin building $X$ of type $(W, S)$.  
\begin{enumerate}[label=(\roman*)]
\item We have $W_{\proj_R(R')} = W_R \cap W_{R'}$. That group is generated by the reflections through the hyperplanes in $\mathcal H(R) \cap \mathcal H(R')$.

\item A chamber $d \in \Ch(R)$ belongs to $\Ch(\proj_R(R'))$ if and only if $d$ and $R'$ are on the same side of every hyperplane $H \in \mathcal H(R) \setminus \mathcal H(R')$. 

\end{enumerate}	
 
\end{lemma}
\begin{proof}
(i) The first assertion follows from \cite[Lemma~3.123]{AbramenkoBrown}, the second from  \cite[Lemma~2.3 and 2.7]{CapRAB}.

\medskip \noindent
(ii) If $H \in \mathcal H(R)$, then some panel $\sigma$ of $H$ is contained in $R$. In view of \cite[2.30.5]{Tits74}, we have $\proj_\sigma(c) = \proj_\sigma(\proj_R(c))$ for all chambers $c$. Thus $c$ and $\proj_R(c)$ are on the same side of $H$. If moreover $H \not \in \mathcal H(R')$, then all chambers of $R'$ are on the same side of $H$. The `only if' part of (ii) follows. 

Let now $d \in \proj_R(R')$ and let $d' \in \Ch(R)$ be any chamber that lies on the same side as $R'$ of every hyperplane $H \in \mathcal H(R) \setminus \mathcal H(R')$. Then the only hyperplanes separating $d$ from $d'$ belong to $\mathcal H(R) \cap \mathcal H(R')$. Since the reflexion through every such hyperplane belongs to $W_{\proj_R(R')}$ by (i), and since that group acts transitively on $\Ch(\proj_R(R'))$, we infer that $d' \in \Ch(\proj_R(R'))$. This proves (ii). 
\end{proof}

Finally, we shall need the following result, that relates the combinatorial projection with the CAT($0$) metric. 

\begin{proposition}\label{prop:proj}
Let $x, x' \in |X|$ be distinct points, let  $R = \Res(x)$, and $ R'=\Res(x')$. The following assertions hold. 
\begin{enumerate}[label=(\roman*)]
\item $[x, x'] \cap |\Ch(R)|$ is a geodesic segment of positive length. 
\item For any $y \in [x, x'] \cap |\Ch(R)|$ with   $x \neq y$, we have 
$$\Ch(x) \cap \Ch(y) = \Ch(\proj_R(R')).$$
\end{enumerate}	
\end{proposition}
\begin{proof}
Since $|\Ch(R)| $ is closed and convex by Lemma~\ref{lem:Ch(R)}, the intersection $[x, x'] \cap |\Ch(R)|$ is a geodesic segment containing $x$. The fact that it has positive length follows from Lemma~\ref{lem:Ch(x)}. This proves (i). 

Let now $y \in [x, x'] \cap |\Ch(R)|$ with $y \neq x$. By definition we have $\Ch(x) = \Ch(R)$.  In order to prove (ii), we need to show that the chamber set of the residue $\proj_R(R')$ coincides with the chamber set of the residue $R \cap \Res(y)$. This is true if and only if that containment holds in some apartment containing $R$ and $R'$. Thus, there is no loss of generality in assuming, for the rest of this proof, that the building $X$ is thin. 

Let us   observe that $\Ch(R) \cap \Ch(y)$ is non-empty. Indeed, since $y \in |\Ch(R)|$, there exists a chamber $d \in \Ch(R)$ with $y \in |d|$. Any such chamber $d$ belongs to the intersection $\Ch(R) \cap \Ch(y)$. We next claim that $d \in \Ch(\proj_R(R'))$.   Let $\mathcal H(R) $ be the collection of all hyperplanes of the thin building $X$ that cross   $\Ch(R)$, and similarly for $R'$. Let $H  \in \mathcal H(R)$ and suppose that $y \in |H|$. Then the entire geodesic segment $[x, x'] $ is contained in the hyperplane  $|H|$ by \cite[Lemma~3.4]{NoskovVinberg}, hence $H \in \mathcal H(R')$ since $R'=\Res(x')$. It follows that, for all hyperplane $H  \in \mathcal H(R) \setminus \mathcal H(R')$, the point $y$ and the residue $R'$ are contained  in the same half-apartment determined by $H$. In particular $H$ does not separate the chamber $d$ from any chamber of $R'$. It follows from Lemma~\ref{lem:proj_R(R')}(i) that $d$ is a chamber of $ \proj_R(R')$, confirming the claim. This shows that $\Ch(R) \cap \Ch(y) \subseteq \Ch(\proj_R(R')).$

Recall from Lemma~\ref{lem:proj_R(R')}(ii) that the stabilizer $W_{T}$ of $T = \proj_R(R')$ is generated by the reflexions through all elements of $ \mathcal H(R) \cap \mathcal H(R')$. Every hyperplane in $ \mathcal H_A(R) \cap \mathcal H_A(R')$ contains $x$ and $x'$, hence also $y \in [x, x']$. Thus $W_T$ fixes $y$. Since $\proj_R(R')$ is a residue of the thin building $X$, it follows that the parabolic subgroup $W_T$ is transitive on the chambers of that residue. Since $W_T$ stabilises $R$, we infer that $ \Ch(R) \cap \Ch(y) \supseteq \Ch(\proj_R(R'))$. The assertion (ii) follows. 
\end{proof}

\subsection{An exhaustion of spherical residues}

Let $X$ be a locally finite building of type $(W,S)$. If $(R_n)_{n\geq 0}$ is a sequence of spherical residues in $X$, denote by $Y_n = \bigcup_{k=0}^n \Ch(R_k)$. We say that the sequence $(Y_n)_{n\geq 0}$ is \textbf{weakly convex} if for all $n\geq 0$ and $b\in Y_n$, we have $\proj_{R_{n+1}}(b) \in Y_n$. We say that it exhausts $X$ if $\bigcup_{n \geq 0} Y_n = X$.

When $(W,S)$ is the free Coxeter system of rank $2$ (i.e. when $X$ is a locally finite tree), then there are many weakly convex sequences in $X$. Indeed, we can argue as follows. First, observe that if $Y\subseteq X$ is any finite subtree, then it has at least one leaf $\ell$. The set of edges in $X$ that contain $\ell$ defines a residue $R$. It is a spherical residue because $X$ is locally finite, and there exists a unique $c\in X$ such that $Y \cap R = \{ c \}$. Now for all $b\in Y$, we have $\proj_R(b) = c \in Y$. Moreover, $Y\cup R$ is also a finite subtree. 

Now we can start with any chamber $c_0 \in X$, and repeatedly use this observation. Set $R_0 = \{ c_0\}$, which is a spherical residue, and $Y_0 = R_0$ is a finite subtree of $X$. By induction, if $Y_n$ is constructed for some $n\geq 0$, choose $R_{n+1}$ as above. The sequence $(Y_n)_{n\geq 0}$ is clearly weakly convex. It is not difficult to picture such a sequence which moreover exhausts $X$, by imposing for instance that the leaf $\ell$ chosen in $Y_n$ is at minimal distance from $c_0$.

When $(W,S)$ is an arbitrary Coxeter system, it is not possible to use the same argument. Indeed, the main difficulty is that unlike the case of trees that we have just discussed, when $Y\subseteq X$ is convex and $R$ is a spherical residue such that $\proj_R(n) \in Y$ for all $b\in Y$, the union $Y\cup R$ need not be convex in general. However, using the interaction between combinatorial projections and the CAT($0$) metric on $|X|$, it is possible to construct a weakly convex exhaustion for arbitrary types.

\begin{proposition}\label{prop:exhaustion}
In every locally finite building  $X$  of type $(W, S)$, there exists a sequence $(R_n)_{n \geq 0}$ of spherical residues satisfying the following conditions, where $  Y_n = \bigcup_{k=0}^n \Ch(R_k)$:
\begin{enumerate}[label=(\arabic*)]
\item $\bigcup_{n \geq 0} \Ch(R_n) = \Ch(X)$.

\item For all $n \geq 0$ and $b \in   Y_n$, we have $\proj_{R_{n+1}}(b) \in Y_n$. 
\end{enumerate}
	
\end{proposition}

\begin{proof}[Proof of Proposition~\ref{prop:exhaustion}]
Pick an arbitrary base point $x_0 \in |X|$. For each chamber $c \in \Ch(X)$, the set $|c|$ is closed and convex. Let $\pi_{|c|} \colon |X| \to |c|$ denote the nearest point projection (see \cite[\S II.2]{BH99}). Let $\mathcal D  = \{d(x_0, \pi_{|c|}(x_0)) \mid c \in \Ch(X)\}$. Since $X$ is locally finite, the intersection of $\mathcal D$ with any finite interval is a finite set. Thus $\mathcal D$ is discrete, and we may enumerate its elements as a strictly increasing sequence, namely $0 = d_0 < d_1<\dots$.  

For each integer $m\geq 0$, we let $\mathcal R_m$ be the set of all spherical residues of the form $\Res(\pi_{|c|}(x_0))$, where $c \in \Ch(X)$ is a chamber such that $d(x_0, \pi_{|c|}(x_0)) = d_m \in \mathcal D$. The required sequence $(R_n)_{n \geq 0}$ is then constructed  by enumerating, successively for all integers $m \geq 0$, the elements of the set $\mathcal R_m$ in an arbitrary order. 

The fact that the sequence $(R_n)_{n \geq 0}$ defined in this way satisfies (1) is clear: for each $c \in \Ch(X)$, the distance $d(x_0, \pi_{|c|}(x_0)) $ is some element $d_m \in \mathcal D$, hence $R = \Res(\pi_{|c|}(x_0))$ belongs to $\mathcal R_m$, so that $c$ belongs to the chamber set of some residue of the sequence $(R_n)$.

It remains to verify the projection condition (2). We shall prove that for all $k \leq n$, the chamber set of the projection $T= \proj_{R_{n+1}}(R_k)$ is entirely contained in $Y_n$. This indeed implies that (2) holds. 

By construction, for each $n$ there is a chamber $c_n \in \Ch(R_n)$ such that $d(x_0, x_n) \in \mathcal D$, where $x_n =  \pi_{|c_{n}|}(x_0)$, and $R_n = \Res(x_n)$. 

Let us now fix $k \leq n$. Set $R = R_{n+1}$, $x=x_{n+1}$, $R' = R_k$, $x'=x_k$ and   $T= \proj_{R}(R')$. Since $k \leq n$, we have $\Ch(R') \subseteq Y_n$. We must show that $\Ch(T) \subseteq Y_n$. 

If $x' = x$ then $R' = R = T$ and there is nothing to prove. We may thus assume that $x \neq x'$. 

Since $k \leq n$, we have $d(x_0, x') = d(x_0, x_k) \leq d(x_0, x_{n+1})=d(x_0, x)$. By the CAT(0) inequality, it follows that $d(x_0, y) < d(x_0, x)$ for all $y \in [x, x']$ different from $x$ and $x'$. By Proposition~\ref{prop:proj}, there exists some $y \in [x, x'] \cap |\Ch(R)|$ different from $x$, and for any such $y$, we have $\Ch(y) \supseteq \Ch(T)$. Thus for each $c \in \Ch(T)$, we have $y \in |c|$, so that $d(x_0, \pi_{|c|}(x_0))\leq d(x_0, y) < d(x_0, x)$. It follows that $c$ belongs to the chamber set of a residue $R_i$ for some $i \leq n$. Hence $c \in Y_n$, as required.
\end{proof}

\section{Hermitian spaces and infinite-dimensional symmetric spaces}
The goal of this section is to describe basic facts on (infinite-dimensional) Hermitian spaces of finite index, and the associated symmetric spaces. The investigation of those objects was initiated by Pontryagin \cite{Pont44}, and stimulated by Gromov~\cite[\S6.A.III]{Gromov}. Our exposition is based on \cite{BIM, DuchesnePhD, Duchesne, MoPy}. 

\subsection{Hermitian spaces of finite index}\label{sec:Hermitian}

Let us review some facts on sesquilinear forms on a  possibly infinite-dimensional $\KK$-vector space $V$, where $\KK \in \{\RR, \CC\}$. We refer to \cite[\S2]{BIM} for more information and background. Let $\varphi \colon V \times V \to \KK$ be a  sesquilinear form. We assume that $\varphi$ is \textbf{symmetric}, i.e. $\varphi(v, w) = \overline{\varphi(w, v)}$ for all $v, w$. The pair $(V, \varphi)$ is then called a \textbf{Hermitian space}. Given a subspace $U$, we write $U^\perp = \{w \in V\mid \varphi(v, w) = 0 \ \forall v \in U\}$. The subspace $V^\perp$ is called the \textbf{radical} of $\varphi$; it is also denoted by $\Rad(\varphi)$. We say that  $\varphi$ is \textbf{non-degenerate} if $\Rad(\varphi) =\{0\}$.

A vector $v \in V$ is \textbf{isotropic} if $\varphi(v, v) = 0$, and \textbf{anisotropic} otherwise. A subspace $W \leq V$ is \textbf{isotropic} if the restriction of $\varphi$ to $W$ vanishes. The supremum of the dimension of isotropic subspaces is called the \textbf{index} of $\varphi$, denoted by $i(\varphi)$.  The \textbf{positive} (resp. \textbf{negative}) \textbf{index} of $\varphi$, denoted by $i_+(\varphi)$ (resp. $i_-(\varphi)$) is the supremum of the dimension of subspaces on which the restriction of $\varphi$ is positive definite (resp. negative definite). 

Suppose that $\varphi$ is non-degenerate. The finite-dimensional theory implies that $i(\varphi)  \leq \min \{i_+(\varphi),i_-(\varphi)\}$. In particular $\varphi$ has finite index as soon as its positive or negative index is finite. In that case we have 
$$i(\varphi) = \min \{i_+(\varphi),i_-(\varphi)\}.$$ 
The ordered pair $(i_-(\varphi), i_+(\varphi))$ is called the \textbf{signature} of $\varphi$.

Still assuming that $\varphi$ is non-degenerate, we denote by $\End(V, \varphi)$ the space of linear endomorphisms $f \in \End(V)$   that have a \textbf{$\varphi$-adjoint} (simply called an \textbf{adjoint} if there is no possible confusion) $f^*$, i.e. a map $f^* \in \End(V)$ such that 
$$\varphi(f(v), w) = \varphi(v, f^*(w))$$
for all $v, w \in V$. The non-degeneracy of $\varphi$ ensures that if $f^*$ exists, it is unique. One may verify that $\End(V, \varphi)$ is a $*$-algebra. 

We shall now define the notion of \textit{strong non-degeneracy}; we will see below that it guarantees the existence of $\varphi$-adjoints for all continuous operators. 

An \textbf{orthogonal decomposition} of $V$ is a direct sum decomposition $V = V_1 \oplus V_2$ such that $\varphi(v_1, v_2) =0$ for all $v_i \in V_i$. In that case we write $V = V_1 \perp V_2$. A \textbf{$\pm$-decomposition} is an orthogonal decomposition $V = V_1 \perp V_2$ such that the restriction of $\varphi$ to $V_1$ is negative definite, and the restriction of $\varphi$ to $V_2$ is positive definite. Such a decomposition gives rise to a positive definite sesquilinear form $\varphi'$ on $V$  whose restriction to $V_1$ (resp. $V_2$) is $- \varphi$ (resp. $ \varphi$), and according to which $V_1$ and $V_2$ remain orthogonal. We say that a non-degenerate form $\varphi$ of finite index (or the Hermitian space $(V, \varphi)$) is \textbf{strongly non-degenerate} if the positive definite Hermitian space $(V, \varphi')$ is a Hilbert space. This does not depend on the choice of a $\pm$-decomposition of $V$; moreover,  all the Hilbertian scalar products respectively associated with those $\pm$-decompositions  are pairwise equivalent , see \cite[Lemma~2.4]{BIM}.

If $\varphi$ is non-degenerate of finite index, there always exists a \textbf{canonical completion} $(\overline V, \overline \varphi)$ which is strongly non-degenerate of finite index, and in which $(V, \varphi)$ embeds as a dense subspace, see \cite[Proposition~3.9]{MoPy2}.
 
 \subsection{Representations on strongly non-degenerate Hermitian spaces}
 
An \textbf{automorphism} of $(V, \varphi)$ is an invertible linear operator $\alpha \colon V \to V$ such that $\varphi(\alpha(x), \alpha(y)) = \varphi(x, y)$ for all $x, y \in V$. The group of all automorphisms is denoted by $\Aut(V, \varphi)$. In case $\varphi$ is strongly non-degenerate of finite index, every automorphism of $(V, \varphi)$ is a bounded linear operator on $V$ with respect to the associated Hilbertian structure (see \cite[\S2]{BIM} and \cite[Proposition~ 9.7]{DuchesnePhD}). 
Given a locally compact group $G$, a linear representation $\sigma \colon G \to \Aut(V, \varphi)$ is called \textbf{continuous} if it satisfies the equivalent conditions highlighted in the following. 

\begin{proposition}\label{prop:continuity}
For a locally compact group $G$ acting by automorphisms on a strongly non-degenerate Hermitian space $(V, \varphi)$ of finite index, the following conditions are equivalent.
\begin{enumerate}[label=(\roman*)]
\item The action map $G \times V \to V$ is continuous. 
\item For each $v \in V$, the orbit map $G \to V \colon g \mapsto gv$ is continuous. 
\item There is a dense subset $W$ in $V$ such that for each $v \in W$, the orbit map $G \to V \colon g \mapsto gv$ is continuous. 
\end{enumerate}
	 
\end{proposition}
\begin{proof}
We have already mentioned that every element of $G$ acts by a continuous operator on $V$. Hence, the equivalence between (i) and (ii) follows from \cite[Proposition 13.2]{Robert_book}. It is clear the (ii) implies (iii). To show the converse, we may assume that $i(\varphi) = i_-(\varphi)$ and consider a $\pm$-decomposition  $V = V_- \perp V_+$. The proof of  \cite[Proposition~ 9.7]{DuchesnePhD} bounds the operator norm of an automorphism $f \in \Aut(V, \varphi)$ by a continuous function of the norm of its restriction of the finite-dimensional subspace $V_-$. It follows that the operator norm of elements of $G$ contained in any fixed compact subset is uniformly bounded. The implication from (iii) to (ii) follows. 
\end{proof}

From now on, the term \textbf{representation} of  a locally compact group $G$ on a strongly non-degenerate Hermitian space of finite index will be used to mean a continuous linear representation. 
	
\medskip
For the rest of this section, we assume that $\varphi$ is 	strongly non-degenerate of finite index. Our next goal is to extend some well known facts from abstract harmonic analysis to this context.  

Observe that, for any $v \in V$, the linear functional $f_v \colon V \to \KK: x \mapsto \varphi(x, v)$ is continuous. Moreover the Riesz representation theorem can be applied separately to each summand in a $\pm$-decomposition of $V$, so that every continuous linear functional on $V$ is of that form. In other words the continuous map $V \to V^* : v \mapsto f_v$ is an isomorphism of $V$ onto its topological dual $V^*$. Similarly, one shows that  every continuous linear operator on $V$ has a $\varphi$-adjoint, see the proof of \cite[Proposition~5.4]{MoPy}.

Let now $\sigma $ be a representation of a locally compact group $G$ on the Hermitian space $(V, \varphi)$. The representation $\sigma$ naturally extends to a $*$-representation of the algebra $C_c^\KK(G)$ of continuous, compactly supported $\KK$ valued functions. Indeed, for $f \in C_c^\KK(G)$, the operator 
$$\sigma(f) = \int_G f(g) \sigma(g)dg$$
is defined as in the classical Hilbert case by requiring that
$$\varphi(v, \sigma(f)w) = \int_G  \varphi(v, f(g)\sigma(g)w)dg$$
for all $v, w \in V$. One checks that this indeed defines a linear representation of $*$-algebras 
$$C_c^\KK(G) \to \End(V, \varphi)$$
that we also denote by $\sigma$. This actually extends to the convolution algebra of compactly supported Radon measures on $G$. In particular, given a compact subgroup $B \leq G$ and $\mu_B$ the normalized Haar measure on $B$, viewed as a Radon measure on $G$, we obtain a continuous idempotent operator $\sigma(\mu_B) $ on $V$.

\begin{proposition}\label{prop:Hermitian-Hecke}
Let $\sigma \colon G \to \Aut(V, \varphi)$ be a representation of a locally compact group $G$ on a strongly non-degenerate Hermitian space of finite index, and let $B \leq G$ be a compact subgroup. 
\begin{enumerate}[label=(\roman*)]
\item The restriction $\varphi_B$ of $\varphi$ on the subspace $V^B$ of $B$-invariant vectors is strongly non-degenerate. 

\item The operator $\sigma(\mu_B)$ acts as the $\varphi$-orthogonal projection on $V^B$. 

\item The space $V^B$ is $C_c^\KK(B\backslash G/B)$-invariant. In particular, the representation $\sigma$ yields a $*$-representation
$\sigma_B \colon C_c^\KK(B\backslash G/B) \to \End(V^B, \varphi_B)$. 
\end{enumerate}
\end{proposition}

\begin{proof}
Observe that   $\mu_B$ is idempotent and self-adjoint; thus the same holds for $\sigma(\mu_B)$. The latter takes its values in $V^B$, and acts as the identity operator on $V^B$. This proves assertion (ii). For (ii), we argue as follows. Suppose that $v \in V^B$ is non-zero. The non-degeneracy of $\varphi$ implies the existence of a vector $w \in V$ with $\varphi(v, w) \neq 0$. Hence $\varphi(\sigma(\mu_B)v, w) = \varphi(v, \sigma(\mu_B)w)$ is non-zero. Hence $\varphi_B$ is non-degenerate. Moreover, $V^B$ is closed in $V$, so in fact we have that $\varphi_B$ is strongly non-degenerate. Assertion (iii) is clear since $C_c^\KK(B\backslash G/B)  = \mu_B *C_c^\KK(B\backslash G/B)*\mu_B$.
\end{proof}
	
A representation on a strongly non-degenerate Hermitian space $(V, \varphi)$ is called \textbf{irreducible}  if the only closed invariant subspaces of $V$ are $\{0\}$ and $V$. 	
The following result plays an important role in relating the irreducible representations of the group $G$ and those of the convolution algebra $C_c^\KK(B\backslash G/B)$. It is well known in the unitary case (see \cite[Ch.~IV, Theorem 2]{Lang}). The finite index Hermitian version, which  is  implicit in  the proof of \cite[Proposition~5.4]{MoPy}, follows from a similar argument as in the unitary case.

\begin{proposition}\label{prop:Hecke-irred}
Let $\sigma \colon G \to \Aut(V, \varphi)$ be a representation of a locally compact group $G$ on a strongly non-degenerate Hermitian space of finite index, and let $B \leq G$ be a compact subgroup. For each closed  $ C_c^\KK(B\backslash G/B)$-submodule $W \leq V^B$, the closure $Z$ of the span of $\sigma(G)W= \bigcup_{w \in W} \sigma(G)w$ is a   $G$-invariant subspace such that $Z^B= W$. 

In particular, if  $\sigma(G)V^B $ spans a dense subspace of $V$, then  $\sigma$ is irreducible if and only if $\sigma_B$ is irreducible. 
\end{proposition}
\begin{proof}
Let $W  \leq V^B$ be a non-zero closed subspace invariant under $ C_c^\KK(B\backslash G/B)$. Let $Z$ be the closure of the span of $\sigma(C^\KK_c(G))W$. Hence $Z$ is a closed $G$-invariant subspace of $V$. The projection $Z\mapsto Z^B$ is given by the operator $ \sigma(\mu_B)$. Hence  $\sigma_B( C_c^\KK(B\backslash G/B)) W= \sigma(\mu_B) \sigma(C^\KK_c(G))W$ spans a dense subspace of $Z^B$. Hence $W = Z^B$. Every closed  $C_c^\KK(B\backslash G/B) $-submodule of $V^B$ is thus attained by the map $Z \mapsto Z^B$. 

Assume now that $\sigma(G)V^B$ spans a dense subspace of $V$. If $W$ is non-zero and if $\sigma$ is irreducible, it follows that the space $Z$ we have just constructed is the whole of $V$. Hence $W = Z^B = V^B$, thereby confirming that $\sigma_B$ is irreducible. 

Assume conversely that $\sigma_B$ is irreducible, and let $Z \leq V$ be a  closed $G$-invariant subspace. Hence $Z^B $ is a closed $\sigma_B$-invariant subspace of $V^B$. The irreducibility yields $Z^B = \{0\}$ or $Z^B = V^B$. 

Suppose first that $Z^B = \{0\}$. By Proposition~\ref{prop:Hermitian-Hecke} the form $\varphi_B$ is strongly non-degenerate, so that $V = V^B \perp (V^B)^\perp$. The equality $Z^B = \{0\}$ yields $Z \leq \ker(\sigma(\mu_B)) = (V^B)^\perp$, so that $V^B \leq Z^\perp$. Since $Z$ is $G$-invariant, we infer that $\sigma(G)V^B  $ is contained in $Z^\perp$, hence $V \leq Z^\perp$ in view of the hypothesis.  Thus $Z \leq V^\perp = \{0\}$. 

Suppose now that $Z^B = V^B$. 
 This implies that the $G$-invariant subspace $Z$ contains $\sigma(G)V^B$, hence $Z= V$ by hypothesis. Hence $\sigma$ is irreducible.
\end{proof}
		
\subsection{Infinite-dimensional symmetric spaces of finite rank}\label{sec:sym-spaces}

Given a strongly non-degenerate Hermitian space $(V, \varphi)$ over $\KK$ of negative index  $p < \infty$ and positive index $q \geq 1$ (possibly $q =\infty$),  the collection $M_{p, q}(\KK)$ of all $p$-dimensional subspaces of $V$ on which the restriction of $\varphi$ is negative definite, is a Riemannian symmetric space, whose dimension is infinite if $\dim(V)$ is infinite. In the latter case $M_{p, q}(\KK)$ is also denoted by $X_p(\KK)$. It carries a natural distance function that turns it to a complete CAT($0$) space, see \cite[Corollary~2.5]{Duchesne}. The group $\Aut(V, \varphi)$ acts continuously by isometries on $M_{p, q}(\KK)$. A systematic study of those spaces in the infinite-dimensional case has been undertaken by B.~Duchesne \cite{DuchesnePhD, Duchesne}. For $p=1$ it is a CAT($-1$) space, called  \textbf{real} or \textbf{complex hyperbolic space} according to whether $\KK= \RR$ or $\CC$.

We need to recall the well known fact that the existence of an unbounded action on a real or complex hyperbolic space of arbitrary dimension is an obstruction to Kazhdan's property (T). We refer to  \cite[Corollary~2.7.3]{BHV} for the finite-dimensional case. The   proof provided by  \cite[Remark~2.11.4]{BHV} actually applies to spaces of arbitrary dimension, see \cite[Theorem~7.4.3]{CCJJV}. 

\begin{proposition}\label{prop:FH}
Let a locally compact group $G$ act continuously by isometries on a   (possibly infinite-dimensional) real  or complex  hyperbolic space. 
\begin{enumerate}[label=(\arabic*)]
\item If $G$ has property (T), then the action has bounded orbits, hence   a fixed point.
\item If the action is metrically proper, then $G$ has the Haagerup property. 
\end{enumerate}
	
\end{proposition}
\begin{proof}
The point is that the distance function on a real or complex hyperbolic space is a kernel conditionally of negative type. Indeed,  this follows from  \cite[Corollary~7.4]{FarHar}   in the finite-dimensional case. The latter implies the result in arbitrary dimension, since the condition concerns finite sets of points, see \cite[Theorem~7.4.3]{CCJJV}.
\end{proof}

\section{Hecke algebras and Hecke-induced modules}\label{sec:Hecke}

\subsection{Iwahori--Hecke algebras}
\label{sec:hecke-preliminaries}

Let $(W, S)$ be a Coxeter system of finite rank and $\KK$ be the field of real or complex numbers. Given an $S$-tuple $\mathbf q = (q_s)_{s \in S}$ of non-zero real numbers such that $q_s = q_t$ if $s$ and $t$ are conjugate, the associated \textbf{Iwahori--Hecke algebra} (or simply the \textbf{Hecke algebra}) is the    $*$-algebra denote by $\KK_{\mathbf q}[W]$ over $\KK$  that has a basis $(e_w)_{w \in W}$ indexed by $W$ and whose multiplication satisfies
$$e_s e_w= \left\{
\begin{array}{ll}
e_{sw} & \text{if } \ell(sw) > \ell(w),\\
q_s e_{sw} + (q_s-1) e_w & \text{if } \ell(sw) < \ell(w)
\end{array}
\right.$$
for all $s \in S$ and $w \in W$, see \cite[\S19.1]{Davis} and \cite[CH.~IV, Exer.~23]{Bour}. The involution $*$ is defined by the assignment 
$$e_w^* = e_{w^{-1}}.$$ 

It is generated by $(e_s)_{s \in S}$. That generating set, together with the following relations, provide a presentation of the Hecke algebra $\KK_{\mathbf q}[W]$:
the \textbf{quadratic relations}
\begin{align*}
 e^2_s & =  (q_s-1) e_s + q_s 
 \end{align*}
 for all $s \in S$, and the \textbf{Artin relations}
 \begin{align*} (e_s e_t)^r& =  (e_t e_s)^r & \text{ if }  m_{s, t} = 2r \\
 (e_s e_t)^r e_s& =  (e_t e_s)^r e_t & \text{ if }  m_{s, t} = 2r+1 
\end{align*}
for all $s, t \in S$ with $m_{s, t}$ finite (see \cite[Ch.~IV, \S 2, Exer.~23]{Bour}).

Let $w \in W$. Write $w$ as a reduced word $w = s_1 \dots s_n$ with $s_i \in S$. We set $q_w = q_{s_1} \dots q_{s_n}$. In view of the solution of the word problem in Coxeter groups, the value of $q_w$ does not depend on the chosen reduced word. It will be convenient to renormalize the bases vectors $e_w$ by setting 
$$T_w = \frac 1 {q_w} e_w.$$
 By \cite[Lemma~19.1.2]{Davis}, we have  $T_u*T_v=T_{uv}$ for all $u, v \in W$ with  $\ell(uv)=\ell(u)+\ell(v)$.  
  In view of \cite[Proposition~19.1.1 and Lemma~19.1.5]{Davis}, we have
\begin{equation}\label{eq:Hecke}
T_w T_s= \left\{
\begin{array}{ll}
T_{ws} & \text{if } \ell(ws) > \ell(w),\\
q_s^{-1} T_{ws} + (1-q_s^{-1}) T_w & \text{if } \ell(ws) < \ell(w)
\end{array}
\right.
\end{equation}
for all $s \in S$ and $w \in W$. In particular, the normalized quadratic relation yields
\begin{equation*}
        T_s^2=q_s^{-1}T_e+(1-q_s^{-1})T_s.
\end{equation*}

\subsection{Hecke algebras and convolution}\label{sec:Hecke-convol}

Let now $X$ be a thick locally finite building of finite rank and type $(W,S)$. From now on, the CAT($0$) metric realization of $X$ will no longer be used and $X$ will be exclusively be viewed as a set of chambers endowed with the $W$-valued  distance $\delta$. For that reason, we will identify $X$ with its chamber set $\Ch(X)$. Similarly, we identify any residue $R$ with $\Ch(R)$. 

We shall assume that $X$ is \textbf{semi-regular}, which means that for each $s \in S$, all $s$-panels of $X$ have the same cardinality. Observe that this condition is automatic if $\Aut(X)$ acts transitively on chambers, since the $\Aut(X)$-action is then transitive on $s$-panels for all $s \in S$.  We   then define $q_s \in \mathbf N$ as the number such that all $s$-panels contain $q_s+1$ chambers. The tuple  $\mathbf q = (q_s)_{s \in S}$ is called the \textbf{thickness vector} of the semi-regular building $X$. 

As in the previous section, we let $\KK$ be the field of real or complex numbers. 
The following lemma recalls the well-known fact that the space $C^\KK_c(X)$ of finitely supported $\KK$-valued functions on $X$ naturally carries the structure of a right  $\KK_{\mathbf q}[W]$-module, see    \cite{Parkinson}. Given a chamber $c \in X$, we denote by $\one_c$ the indicator function $\one_{\{c\}} \colon X \to \KK$. We do not use the symbol $\delta$, since it is already used to denote the Weyl distance on $X$. 
For $\KK = \CC$, the exponent $\KK$ is often omitted when this does not cause any confusion. 

\begin{lemma}\label{lem:right-H-mod}
Let $\qq = (q_s)_{s \in S}$ be the thickness vector of $X$. 
  Given $x \in X$ and $s \in S$, we set  
$$\one_x \cdot e_s =   \sum_{y \in X, \delta(x, y)=s} \one_y. $$
Those assignments endow $C^\KK_c(X)$ with the structure of a right $\KK_{\mathbf q}[W]$-module. 
\end{lemma}
\begin{proof}
Since the functions $(\one_x)_{x \in X}$ form a basis of  $C_c^\KK(X)$, the given assignments define a right action of the free associative $\KK$-algebra generated by $(e_s)_{s \in S}$ on the space $C_c^\KK(X)$. We need to verify that the definition relations of $\KK_\qq[W]$ are satisfied.
We observe that 
\begin{align*}
(\one_x \cdot e_s)\cdot e_s & = \sum_{y \in X, \delta(x, y)=s} \one_y \cdot e_s\\
& =  \sum_{y \in X, \delta(x, y)=s} \sum_{z \in X, \delta(y, z) = s} \one_z \\
&= q_s \one_x + (q_s-1) \sum_{z \in X, \delta(x, z)=s} \one_z\\
&=q_s \one_x + (q_s-1) \one_x \cdot e_s\\
&= \one_x \cdot e_s^2,
\end{align*}
hence the quadratic relations are satisfied. 

For $s \neq t \in S$ with $m_{s, t}$ finite, let $w_0 $ be the longest word in $\langle s, t \rangle$. Let also $f_{s, t} = e_s e_t \dots$ be the product of $m_{s, t}$ generators starting with $e_s$ on the left, and define $f_{t, s}= e_t e_s \dots$ similarly. It follows from the definitions that, for all $x \in X$, we have
 \begin{align*}
\one_x \cdot  f_{s, t} &=  \sum_{y \in X, \delta(x, y)=w_0} \one_y 
= \one_x \cdot f_{t, s},
\end{align*}
hence the Artin relations are also satisfied. The required assertion follows. 
\end{proof}

\medskip
Let now $G$ be a locally compact group acting properly, continuously and transitively by automorphisms on $X$.
Fix a base chamber $c_0$. The transitivity of $G$ implies that  it is generated by the stabilizers of the panels containing $c_0$, which are all compact and open. In particular $G$ is unimodular.

Let $B=\Stab_G(c_0)$, which  is also a compact open subgroup of $G$. We normalize Haar measure $m$ on $G$ by imposing that $m(B)=1$. The orbit map defines a one-to-one correspondence between the coset space $G/B$ and $X$. Therefore the space $C^\KK_c(X)$ can be identified with the subspace $C^\KK_c(G/B)$ of $C^\KK_c(G)$ consisting of the right-$B$-invariant functions. We view $C^\KK_c(G)$, equipped with the convolution product,  as a $\KK$-algebra  endowed with the standard involution  $f \mapsto f^*$ defined by 
$$f^*(x) = \bar{f}(x^{-1}).$$ 
The subalgebra consisting of bi-$B$-invariant functions is denoted by $C^\KK_c(B\backslash G/B)$. For each $w \in W$, we define
$$D(w) = \{g \in G \mid \delta(c_0, gc_0) = w\}.$$
Since $B$ fixes $c_0$, it follows that $D(w)$ is a union of double $B$-cosets. In other words the characteristic function $\one_{D(w)}$ is bi-$B$-invariant. The building $X$ being locally finite, we have $\one_{D(w)} \in C^\KK_c(B\backslash G/B)$. Moreover, we have
$$m(D(w)) = q_w$$
for all $w \in W$, since the Haar measure $m$ has been normalized in such a way that $m(B)=1$. 

By \cite[Lemma~19.1.5]{Davis}, the map 
$$\KK_\qq[W] \to C^\KK_c(B\backslash G/B) : e_w \mapsto \one_{D(w)}$$
sending $e_w$ to the characteristic function $\one_{D(w)} \in C^\KK_c(G)$ establishes an isomorphism of the Hecke algebra $\KK_{\mathbf q}[W]$ onto a subalgebra $J$ of $C^\KK_c(B\backslash G/B)$.  We say that the $G$-action on $X$ is \textbf{Weyl-transitive} if, for each $w\in W$, the set $D(w)$ is a single $B$-double coset, namely $BwB$. In that case we have $m(BwB) = q_w$ for all $w \in W$, and the map  $\KK_\qq[W] \to C^\KK_c(B\backslash G/B)$  is surjective, so that $\KK_\qq[W] $ and $ C^\KK_c(B\backslash G/B)$ are isomorphic.

The convolution product naturally endows the $\KK$-vector space $C^\KK_c(G/B)$ with the structure of a left $C^\KK_c(G)$-module and a right $C^\KK_c(B\backslash G/B)$-module. 
Since $G/B$ can be identified with $X$, the spaces $C^\KK_c(G/B)$ and $C^\KK_c(X)$ are naturally isomorphic. Viewing   $\KK_{\mathbf q}[W]$ as the subalgebra $J$ of $C^\KK_c(B\backslash G/B)$ as above, this recovers  Lemma~\ref{lem:right-H-mod}  in this case.

\subsection{Hermitian modules and Hecke-induction}\label{sec:Hecke-ind}

Let   $(V_0, \varphi_0)$ be a non-degenerate Hermitian space over $\KK$,  and 
$$\rho \colon \KK_{\mathbf q}[W]\to \End(V_0, \varphi_0)$$ 
be a linear representation of $*$-algebras. Thus $\rho$ is a homomorphism of algebras, such that $\rho(a^*) = \rho(a)^*$ for all $a \in  \KK_{\mathbf q}[W]$. We say that $(V_0, \varphi_0)$ is a \textbf{Hermitian $\KK_{\mathbf q}[W]$-module}.  The action map $\rho$ will often be omitted.

 Let $X$ be a   locally finite building of type $(W, S)$ and suppose that $\qq$ is its thickness vector. In that case, we define the \textbf{Hecke-induced module} associated with the Hermitian $\KK_\qq[W]$-module $(V_0, \varphi_0)$ as the tensor product of modules
$$\Ind(V_0) = C^\KK_c(X) \otimes_{\KK_\qq[W]} V_0,$$
where $\qq=(q_s)_{s \in S}$ is the thickness vector of $X$ (see Lemma~\ref{lem:right-H-mod} for the right $\KK_\qq[W]$-module structure on $C^\KK_c(X)$). To define a linear left action of $\Aut(X)$, it is convenient to consider the tensor product of vector spaces 
$$\mathscr E = C^\KK_c(X) \otimes V_0.$$
The induced module is thus a quotient of $\mathscr E$. The full automorphism group $\Aut(X)$ acts by linear transformations on $\mathscr E$: the element $g \in \Aut(X)$ maps $\one_c \otimes v$ to $\one_{g(c)} \otimes v$ for all $c \in \Ch(W)$ and $v \in V_0$. It is straightforward to check that this action preserves the kernel of the quotient map $\mathscr E \to \Ind(V_0)$. Hence it defines an $\Aut(X)$-action by linear transformations on the induced module $\Ind(V_0)$. 

The   form $\varphi_0$ naturally induces a Hermitian form on the induced module. In order to define it, we first consider the sesquilinear form $\varphi$  on $\mathscr E$ defined by 
$$\varphi(\one_{c_1} \otimes v_1, \one_{c_2} \otimes v_2) = \varphi_0(T_{\delta(c_2, c_1)}v_1, v_2)$$
for all $c_1, c_2 \in \Ch(X)$ and $v_1, v_2 \in V_0$. 
The following result ensures that $\varphi$ canonically descends to a Hermitian form on $\Ind(V_0)$. We denote it by $\Ind(\varphi_0)$, or by $\varphi$ when it does not cause confusion, and call it  the \textbf{induced form}.
Since $\Aut(X)$ preserves the Weyl-distance, that action preserves the form $\varphi$ on $\mathscr E$, hence the induced form on the induced module $\Ind(V_0)$.

\begin{lemma}\label{lem:induced-form}
For all $x \in X$, $v \in V_0$  and $a \in \KK_\qq[W]$, the element 
$$\one_x  \otimes a v - (\one_x \cdot a)\otimes  v  $$ 
belongs to the radical of $\varphi$. 
\end{lemma}
\begin{proof}
It suffices to verify the claim for $a=T_s$ with $s \in S$. Let $y \in X$ and $w \in V_0$. We have
$$
\varphi(\one_x \otimes T_s v, \one_y\otimes w) = \varphi_0(T_{\delta(y, x)} T_s v, w).$$
In view of (\ref{eq:Hecke}), 
the latter is equal to 
$$ \varphi_0(T_{\delta(y, x) s} v, w) $$
 if $\ell(\delta(y, x)s) > \ell(\delta(y, x)) $, and to 
$$ q_s^{-1} \varphi_0(T_{\delta(y, x) s} v, w)+(1-q_s^{-1})\varphi_0(T_{\delta(y, x)}v, w)$$
if   $\ell(\delta(y, x)s) < \ell(\delta(y, x))$. 

Let $x' $ be the projection of the chamber $y$ on the $s$-panel containing $x$, that we denote by $\sigma$. Moreover, we have
\begin{align*}
\varphi((\one_x \cdot T_s)\otimes   v, \one_y\otimes w) & = q_s^{-1}\sum_{z \in \sigma, z \neq x}\varphi(\one_z \otimes v, \one_y\otimes w)\\
&= q_s^{-1}\sum_{z \in \sigma, z \neq x}\varphi_0(T_{\delta(y, z)}v, w).
\end{align*}
Since $\sigma$ contains $q_s+ 1$ chambers, the right-hand side  is equal to 	
$$ \varphi_0(T_{\delta(y, x) s} v, w) $$
if $x' = x$, and to 
$$ q_s^{-1} \varphi_0(T_{\delta(y, x) s} v, w)+q_s^{-1}(q_s-1)\varphi_0(T_{\delta(y, x)}v, w)$$
if $x' \neq x$. 

Since $x' =x$ if and only if $\ell(\delta(y, x)s) > \ell(\delta(y, x)) $, we deduce  that 
$$\varphi(\one_x \otimes T_s v, \one_y\otimes w) = \varphi((\one_x \cdot T_s)\otimes   v, \one_y\otimes w)$$ 
for all $x, y \in X$, $s \in S$ and $v, w \in V_0$. The lemma follows. 
\end{proof}
		
The following basic observation records when the Hecke-induced module is a multiple of the trivial representation. 

\begin{lemma}\label{lem:trivial-induced}
Let   $(V_0, \varphi_0)$ be a non-degenerate Hermitian $\KK_\qq[W]$-module, and let $X$ be a locally finite, semi-regular, building of type $(W, S)$ whose thickness vector is $\qq$.  If $T_s$ acts as the identity operator on $V_0$ for all $s \in S$,  the induced module $(\Ind(V_0), \Ind(\varphi_0))$ is  endowed with the trivial $\Aut(X)$-action.  Moreover the quotient $\Ind(V_0)/\Rad(\Ind(\varphi_0))$ is isomorphic to $V_0$.
\end{lemma}
\begin{proof}
By hypothesis we have $\varphi(\one_x \otimes v , \one_y \otimes w) = \varphi_0(v, w)$ for all $x, y \in X$ and $v, w \in V_0$. Choose a base chamber $c \in X$. It follows that the vector $(\one_c -   \one_x) \otimes v$ belongs to $\Rad(\varphi)$ for all $x \in X$. This implies that the restriction of the quotient map $\mathscr E \to \Ind(V_0)/\Rad(\Ind(\varphi_0))$ to the subspace $\one_c \otimes V_0$ is surjective. The hypothesis that $\varphi_0$ is non-degenerate implies that it is also injective. All the required assertions follow. 
\end{proof}

In the special case where $X$ is locally finite and $\Aut(X)$ contains a closed subgroup $G$ that is transitive on $X$, the Hecke-induced module has another natural  realization. Indeed, as recalled in \S\ref{sec:Hecke-convol}, we have isomorphism $\KK_{\mathbf q}[W] \to J$  is defined by the assignments $e_w \mapsto \one_{D(w)}$ for all $w \in W$.  We define a map 
$$P \colon C^\KK_c(G) \to J : f \mapsto  \sum_{w \in W} \frac 1 {q_w}\langle f, \one_{D(w)} \rangle \one_{D(w)}.$$ 
Thus the value of $f$ at $g \in D(w)$ is equal to the average of the values of $f$ on $D(w)$. 

\begin{lemma}\label{lem:induced-transitive-case}
Let   $X$ be a locally finite, semi-regular, building of type $(W, S)$ whose thickness vector is $\qq$, and $G \leq \Aut(X)$ be a closed subgroup acting transitively. Let also $(V_0, \varphi_0)$ be a   Hermitian $\KK_\qq[W]$-module. 
\begin{enumerate}[label=(\roman*)]
\item The induced module $\Ind(V)$ is isomorphic to $C_c^\KK(G) \otimes_J V_0$. 
\item The induced form $\Ind(\varphi_0)$ coincides with the form on  $C_c^\KK(G) \otimes_J V_0$ defined by 
$$\varphi(f_1 \otimes v_1, f_2 \otimes v_2) = \varphi_0(P(f_2^* * f_1)v_1, v_2).$$
\item If $G$ is Weyl-transitive, then $J = C^\KK_c(B\backslash G/B)$ and 
we have 
$$P(f) = \one_B * f * \one_B$$ 
for all $f \in C_c^\KK(G)$.  
\end{enumerate}
\end{lemma}
\begin{proof}
Let us identify $\KK_\qq[W]$ with $J$ via the isomorphism recalled above. In particular $T_e = \one_{D(e)} = \one_B$. 

Since $G$ is chamber-transitive, the left $G$-modules $C_c^\KK(X)$ and $C_c^\KK(G/B)$ are naturally isomorphic. Moreover, the natural isomorphism also preserves the structure of a right $\KK_\qq[W]$-module. Assertion~(i) follows. 

In order to prove the assertion~(ii),  it suffices to check the equality for characteristic functions of left $B$-cosets, say $f_i = \one_{g_iB}$ with $g_1, g_2 \in G$, since those form a basis of $C_c^\KK(G/B)$. We have $f_2^* * f_1 =m(Bg_2^{-1}g_1B)^{-1} \one_{Bg_2^{-1}g_1B}$. Let $w = \delta(c_0, g_2^{-1} g_1 (c_0)) = \delta(g_2(c_0), g_1(c_0))$. Since $Bg_2^{-1}g_1B \subseteq D(w)$, we  have
\begin{align*}
P(f_2^* * f_1) &= \frac 1 {q_w} \frac 1 {m(Bg_2^{-1}g_1B)} \big(\int_{D(w)} \one_{Bg_2^{-1}g_1B}(x) dx\big) \one_{D(w)}\\
&= \frac 1 {q_w} \one_{D(w)}\\
&= T_w.
\end{align*}
Therefore we obtain that $\varphi_0(P(f_2^* * f_1)v_1, v_2) = \varphi_0(T_{\delta(g_2(c_0), g_1(c_0))}v_1, v_2)$, which confirms (ii). 

For (iii), observe that for all $g \in D(w)$, we have 
$$\one_B * f * \one_B(g) = \int_B \int_B f(xgy)dx dy =\frac 1 {m(BwB)} \int_{BwB} f(h)dh.$$ 
The assertion follows since we have $BwB = D(w)$ if $G$ is Weyl-transitive. 
\end{proof}

\subsection{Negative index of the induced form}
We shall assume henceforth that $(V_0, \varphi_0)$ is a non-degenerate Hermitian $\KK_\qq[W]$-module, and that the negative index 
$$ p = i_-(\varphi_0)$$
is finite.  In general, the induced form $\Ind(\varphi_0)$  can be degenerate. 
 Our goal is to find conditions ensuring that its negative index is also equal to $p$. 

We retain the notation of \S\ref{sec:Hecke-ind}. 
Given a subset $Y \subseteq X$, we further define $\mathscr E_Y$ as the natural image of $C^\KK_c(Y) \otimes V_0$ in $\mathscr E$, and view it as a subspace of $\mathscr E$. The restriction of $\varphi$ to  $\mathscr E_Y \times \mathscr E_Y$ is denoted by $\varphi_Y$.

\begin{proposition}\label{prop:add-residue}
Let $Y \subseteq X$ be a non-empty subset and $R$ be a residue such for all $y \in Y$, we have 
$$\proj_R(y) \in Y.$$
If $i_-(\varphi_0) = i_-(\varphi_R) = i_-(\varphi_Y) = p< \infty$, then $i_-(\varphi_{R \cup Y}) = p$.
\end{proposition}
\begin{proof}
We consider the linear operator $L$ on $\mathscr E_{R \cup Y}$ defined by the assignments
$$\one_c \otimes x \mapsto \left\{
\begin{array}{ll}
\one_c \otimes x & \text{if } c \in R\\
\one_c  \otimes x - \one_{c'} \otimes  T_{\delta(c', c)}x & \text{if } c \in Y \setminus R,
\end{array}
\right.
$$
for all $c \in R \cup Y$ and $x \in V_0$, where $c' = \proj_R(c)$. Observe that $L$ is invertible. Moreover it acts  as the identity operator on $\mathscr E_R$. We set $\mathscr F = L(\mathscr E_{Y \setminus R})$. Since $\mathscr E_{R \cup Y} \cong \mathscr E_{R} \oplus \mathscr E_{Y \setminus R} $, we have $\mathscr E_{R \cup Y} \cong \mathscr E_{R} \oplus \mathscr F$.

We claim that $\mathscr E_R \perp \mathscr F$. Indeed, given $c \in Y \setminus R$, $d \in R$ and $x, y \in V_0$, we set $c' = \proj_R(c)$ and we obtain 
\begin{align*}
\varphi(L(\one_c \otimes x), \one_d \otimes y) 
&= \varphi(\one_c \otimes x -  \one_{c'} \otimes T_{\delta(c', c)}x, \one_d \otimes y) \\
&= \varphi_0(T_{\delta(d, c)}x, y) - \varphi_0(T_{\delta(d, c')} T_{\delta(c', c)}x, y)\\
&= \varphi_0(T_{\delta(d, c)}x, y) - \varphi_0(T_{\delta(d, c')} T_{\delta(c', c)}x, y)\\
&= \varphi_0(T_{\delta(d, c)}x, y) - \varphi_0(T_{\delta(d,   c)}x, y)\\
&=0,
\end{align*}
because $\dist(d, c) = \dist(d, c') +\dist(c', c)$ since 	 $c' = \proj_R(c)$. This proves the claim. 

We have just obtained an orthogonal decomposition  $\mathscr E_{R \cup Y} = \mathscr E_R \perp  \mathscr F$.  The hypotheses imply that $R \cap Y$ contains all projections $\proj_R(y)$ with $y \in Y$. This implies that $\mathscr F$ is a subspace of $\mathscr E_Y$, hence we infer that $\mathscr E_Y = \mathscr E_{R \cap Y}  \perp \mathscr F$. Moreover, it also implies that $R \cap Y$ is non-empty. Hence it contains a subspace $\mathscr E_{\{c\}}$ with $c \in R \cap Y$. That subspace, endowed with the restriction of $\varphi$,  is isomorphic to $(V_0, \varphi_0)$. Hence $p = i_-(\varphi_0) \leq i_-(\varphi_{R \cap Y})$. By hypothesis we have
$$p = i_-(\varphi_Y) =   i_-(\varphi_{R \cap Y}) + i_-(\varphi|_{\mathscr F \times \mathscr F}) \geq p + i_-(\varphi|_{\mathscr F \times \mathscr F}). $$
It follows that the negative index of $\varphi|_{\mathscr F \times \mathscr F}$ vanishes. In other words $\varphi|_{\mathscr F \times \mathscr F}$  is positive semi-definite.

Invoking again the claim, the hypothesis that $i_-(\varphi_R)=p$ now yields
$$i_-(\varphi_{R \cup Y}) = i_-(\varphi_R) + i_-(\varphi|_{\mathscr F \times \mathscr F}) = p +0 = p, $$ 
as required. 
\end{proof}
		
The following subsidiary fact will allow us to check that the condition on $\varphi_R$ in Proposition~\ref{prop:add-residue} is satisfied by spherical residues. 
		
\begin{lemma}\label{lem:spherical-residue-positivity}
Let $X$ be a locally finite building of \emph{spherical} type $(W, S)$, whose automorphism group is transitive, and let $\mathbf q$ be its thickness vector. Let also  $\rho \colon \KK_{\mathbf q}[W]\to \End(V_0, \varphi_0)$ be a linear representation of $*$-algebras. If the form $\varphi_0$ is positive definite, then the induced form $\Ind(\varphi_0)$ is positive semi-definite. 
\end{lemma}

For the proof, we shall rely on \cite[Lemma~1.7]{Rieffel} ensuring the positive semi-definiteness of  the induced form that equips  induced modules of $C^*$-algebras. As in \cite[Definition~1.3]{Rieffel}, given a $C^*$-algebra $A$ and a $C^*$-subalgebra $J$, we say that a map $P \colon A \to J$ is a \textbf{conditional expectation} if $P$ is a continuous, positive, projection such that $P(ab) = P(a)b$ for all $a \in A$ and $b \in J$. A typical example is afforded by the map   $P(a) = pap$, where $p \in A$ is a self-adjoint idempotent and $J = pAp$.  

\begin{lemma}\label{lem:f-d-hilbert-algebra}
Let $A$ be a finite-dimensional $C^*$-algebra and $J$ be a $C^*$-subalgebra. Suppose that $A$ is equipped with an inner product $\langle \cdot, \cdot \rangle$ that makes it a Hilbert algebra (see \cite[A.54]{Dixmier}). Then the orthogonal projection $A \to J$ is a conditional expectation. 
\end{lemma}
\begin{proof}
First observe that for all $x \in J$ and $y \in J^\perp$, we have $xy \in J^\perp$: indeed, for all $z \in J$ we have
$$\langle z, xy \rangle = \langle x^*z, y \rangle =0$$
since $x^*z \in J$. Next observe that  $J^\perp$ is $*$-invariant, since $\langle x, y \rangle = \langle y^*, x^* \rangle$ for all $x, y \in A$ by the definition of a Hilbert algebra. We deduce that for all $x \in J$ and $y \in J^\perp$, we also have  $yx \in J^\perp$. 

Let now $P \colon A \to J$ be the orthogonal projection. Given  $a \in A$, let $a = a_0 + a_1$ be the decomposition of $a$ according to $A = J \oplus J^\perp$. For all $b \in J$, we have 
$$P(ab) = P(a_0 b) + P(a_1 b) = P(a_0 b) = a_0 b = P(a)b$$
since $a_1b \in J^\perp$ be the initial observation made above. 

It remains to prove that $P$ is positive. Since $J$ and $J^\perp$ are both $*$-invariant, we have $P(a)^* = P(a^*)$ for all $a \in A$. Let us now assume that $a \in A$ is  positive. Hence   $P(a)$ is self-adjoint, and in order to prove   that $P(a)$ is positive, it suffices to show that the eigenvalues of the operator on the definite-dimensional algebra $J$ defined by left multiplication by $P(a)$ are non-negative. Let thus $
x \in J\setminus \{0\} $ and $\lambda \in \RR$ with $P(a)x = \lambda x$. We have seen that $P(ax) = P(a)x$. Setting $y = ax - P(ax) \in J^\perp$ and using that $a$ is positive, we observe that  
$$0 \leq \langle ax, x \rangle = \langle P(ax) + y,  x \rangle = \langle \lambda x+y, x\rangle = \lambda \langle   x, x\rangle.$$
This confirms that $\lambda \geq 0$, thereby proving that $P(a)$ is positive.  
\end{proof}

An alternative, more conceptual, proof of the lemma can be derived from the fact that the Hilbert algebra inner product is induced by a faithful trace: indeed, the assertion   follows from the standard realization of
trace-preserving conditional expectations as orthogonal projections \cite[Theorem~9.1.2 and Remark~9.1.3]{AnantharamanPopa}.

\begin{proof}[Proof of Lemma~\ref{lem:spherical-residue-positivity}] 
The chamber-transitivity implies that $\KK_{\mathbf q}[W]$ is isomorphic to a subalgebra $J$ of $C_c(G)$: as recalled in \S\ref{sec:Hecke-convol}, the isomorphism $\KK_{\mathbf q}[W] \to J$ is defined by the assignments $e_w \mapsto \one_{D(w)}$ for all $w \in W$. 

Since $G$ is finite, the algebra $C_c(G) = \ell^2(G)$ is a Hilbert algebra with respect to the form $\langle f, g\rangle = \int_G f(x)\bar g(x) dx$   (see \cite[\S13.10.1]{Dixmier}). Moreover $J$ is a subalgebra admitting an orthonormal basis afforded by the vectors 
$$U_w = q_w^{-\frac 1 2} \one_{D(w)}$$
for all $w \in W$ (see \cite[\S19.2]{Davis}). 

Observe that the map $P \colon C_c(G) \to J : f \mapsto \frac 1 {q_w}\langle f, \one_{D(w)} \rangle \one_{D(w)}$ defined in \S\ref{sec:Hecke-ind} coincides with 
$$P \colon C_c(G) \to J : f \mapsto \sum_{w \in W} \langle f, U_w\rangle U_w. $$
In particular $P$ is the orthogonal projection to $J$ with respect to the form $\langle \cdot, \cdot \rangle$. By Lemma~\ref{lem:f-d-hilbert-algebra}, the map $P$ is a conditional expectation.  Hence the conclusion of the lemma over complex numbers follows from  \cite[Lemma~1.7]{Rieffel}.   The statement over the reals follows by restricting scalars. 
\end{proof}

We are now ready to prove the main result, stated as Theorem~\ref{thmintro:signature} in the introduction. For the reader's convenience, we repeat its statement here.  

\begin{theorem}\label{thm:main}
Let $X$ be a thick, locally finite, semi-regular building of type $(W, S)$, all of whose spherical residues have a  transitive automorphism group.  Let $\mathbf q = (q_s)_{s\in S}$ be its thickness vector, and let $(V_0, \varphi_0)$ be a non-degenerate Hermitian $\KK_{\mathbf q}[W]$-module. 

Suppose that $i_-(\varphi_0) = p<\infty$. If for each spherical subset $J  \subset S$, there is a   $p$-dimensional subspace $U_J$ of $V_0$ that is pointwise fixed by $T_s =  q_s^{-1} e_s$  for all $s \in J$, and on which the restriction of $\varphi_0$ is negative definite, then   $i_-(\Ind(\varphi_0)) = p$. 
\end{theorem}
\begin{proof} 
We  first claim that, for each spherical residue $R$, the form $\varphi_R$ on $\mathscr E_R$ has negative index $p$. 

Let $J \subseteq S$ be the type of $R$, and $H_J$ be the subalgebra of $\KK_{\mathbf q}[W]$ generated by $\{T_s \mid s \in J\}$. Hence $V_0$ is a Hermitian $H_J$-module, and $\mathscr E_R$ is the associated induced module. We have  an orthogonal sum decomposition $V_0 = U_J \perp  U'_J$ of Hermitian $H_J$-modules, where $U_J$ is the $p$-dimensional subspace afforded by the hypotheses. We know   that the restriction of $\varphi_0$ on $U_J$ is negative definite, hence its restriction on $U'_J$ is positive definite.  

The $H_J$-invariant decomposition  $V_0= U_J \perp  U'_J$ yields  an orthogonal sum decomposition $\mathscr E_R = \mathscr E_J \perp \mathscr E'_J$ of the space    $\mathscr E_R = C_c(R) \otimes V_0$, that is preserved by the natural action of $\Aut(R)$.  We denote by $\varphi_R$, $\varphi_J$ and $\varphi'_J$ the restrictions of the form $\varphi$ on the subspaces  $\mathscr E_R$, $\mathscr E_J $ and $ \mathscr E'_J$.  

 Since every generator $T_s$ of $H_J$ acts by the identity operator on $U_J$, it follows from Lemma~\ref{lem:trivial-induced} that the quotient     $\mathscr E_J/\Rad(\varphi_J)$ is   isomorphic to $U_J$, and that the induced   $\Aut(R)$-action is trivial.  In particular we have $i_-(\varphi_J) = i_-(\varphi_0) = p$.

Let us now consider the $H_J$-module $U'_J$. We denote by $\Ind(U'_J)$ the Hecke-induced module for the group $\Aut(R)$. By Lemma~\ref{lem:induced-form}, the kernel of quotient map $ \mathscr E'_J \to \Ind(U'_J)$ is contained in the radical of the restriction $\varphi'_J$ of $\varphi_R$ to $\mathscr E'_J$. In particular the negative index of  $ \varphi'_J$ coincides with the negative index of the induced form on $\Ind(U'_J)$.  Since $\varphi_0$ restricts to a positive definite form on $U'_J$, the induced form on $\Ind(U'_J)$ is positive semi-definite by Lemma~\ref{lem:spherical-residue-positivity}. Hence $\varphi'_J$ is positive semi-definite. Thus $i_-(\varphi_R) = i_-(\varphi_J) + i_-(\varphi'_J) = p +0= p$, which confirms the claim.

\medskip
We now consider the sequence $(R_n)_{n \geq 0}$ of spherical residues afforded by Proposition~\ref{prop:exhaustion}. For each $n \geq 0$ we set $Y_n = \bigcup_{k =0}^n R_k$. Since $Y_0 = R_0$ is a spherical residue, the restriction $\varphi_{Y_0}$ has negative index $p$ by the claim. It now follows by induction on $n$ that the restriction $\varphi_{Y_n}$ has negative index $p$ for all $n$: indeed, the induction step is afforded by Proposition~\ref{prop:add-residue}, whose hypotheses are satisfied in view of Proposition~\ref{prop:exhaustion} combined with the claim above. The conclusion follows. 
\end{proof}

Specializing to $p=0$, we obtain the following consequence, which generalizes Theorem~\ref{thmintro:unitary} from the introduction. 

\begin{corollary}\label{cor:main}
Let $X$, $(W, S)$ and $\mathbf q$ be as in Theorem~\ref{thm:main}. Let $(V_0, \varphi_0)$ be a Hermitian $\KK_{\mathbf q}[W]$-module. 

If $\varphi_0$ is positive definite, then the induced form $\varphi$ is positive semi-definite. 
\end{corollary}

\begin{proof}[Proof of Corollaries~\ref{cor:equivalence} and~\ref{corintro:corner}]
First observe that the map $\mathcal H \mapsto \mathcal H^B$ preserves irreducibility of representations by Proposition~\ref{prop:Hecke-irred}. 

Let us show that every irreducible unitary Hecke-module is attained. Indeed, given a Hermitian $\CC_{\mathbf q}[W]$-module $(V_0, \varphi_0)$ with $\varphi_0$ positive definite, we obtain a unitary representation $\sigma$ of $G$ on the canonical completion $\mathcal V$ of the induced module $\Ind(V_0)$ divided by the radical of $\Ind(\varphi_0)$. It is continuous by Proposition~\ref{prop:continuity}, since $\mathcal V$ contains a dense set of vectors whose stabilizer in $G$ is open. It follows from Proposition~\ref{prop:Hecke-irred} that if  $V_0$ is irreducible, then $\sigma$ is irreducible and the associated representation of $\CC_{\mathbf q}[W] \cong C_c(B\backslash G/B)$ on $\mathcal V^B$ is isomorphic to $V_0$. The required conclusions follow.
\end{proof}

\subsection{The geometric Hecke module}\label{sec:geom-Hecke-module}
					
Given a Coxeter system $(W, S)$, the group $W$ has a faithful linear representation on a real vector space $\bigoplus_{s \in S} \RR \alpha_s$ of dimension $|S|$ preserving a symmetric bilinear form $\varphi_{\mathrm T}$, called the \textbf{Tits form}, defined by 
$$\varphi_{\mathrm T}(\alpha_s, \alpha_t) = -\cos(\frac \pi {m_{s, t}})$$
for all $s, t \in S$, 
where $m_{s, t}$ is the order of $st$, see \cite[Ch.~V, \S 4]{Bour}. That representation is called the \textbf{geometric representation} of $W$. 

As before, we denote by $\KK$ the field of real or complex numbers, and let $\qq$ be an $S$-tuple of elements of $\KK$. The Hecke algebra $\KK_{\mathbf q}[W]$ can be viewed as a deformation of the group algebra $\KK[W]$, which is isomorphic to the Hecke algebra associated with the parameters $q_s = 1$ for all $s$. For an arbitrary $\qq$, the  algebra $\KK_{\mathbf q}[W]$  also admits linear representations on the $\KK$-vector space $V_{\mathrm g} = \bigoplus_{s \in S} \KK \alpha_s$. In order to define it, we 
shall assume that $m_{s, t} \in \{2, 3, 4, 6, 8, \infty\}$ for all $s \neq t$. The latter condition is automatic if $W$ is the Weyl group of a thick locally finite building, in view of the well known restrictions on the arity of the finite generalized polygons, see \cite{FeitHigman}.  

We assume that the tuple $\mathbf q = (q_s)_{s \in S}$ consists of positive reals. For all $s, t \in S$ with $m_{s, t} = \infty$, we make a choice of a scalar $\mu_{s, t} \in \KK$ in such a way that $\mu_{t, s} = \overline{\mu_{s, t}}$. 

For all $s, t \in S$, we set
\begin{equation}\label{eq:eps_st} 
\varepsilon_{s, t} = \left \{
\begin{array}{ll}
-1 & \text{if } s = t\\
0 & \text{if } m_{s, t} = 2\\
\frac{\sqrt{q_s}}{1+q_s} =  \frac{\sqrt{q_t}}{1+q_t} & \text{if } m_{s, t} = 3 \\
\sqrt{\frac{ q_s+q_t}{(1+q_s)(1+q_t)}} & \text{if } m_{s, t} = 4\\
\sqrt{\frac{ q_s+q_t+\sqrt{q_s q_t}}{(1+q_s)(1+q_t)}} & \text{if } m_{s, t} = 6\\
\sqrt{\frac{ q_s+q_t+\sqrt{2q_s q_t}}{(1+q_s)(1+q_t)}} & \text{if } m_{s, t} = 8\\
\mu_{s, t} & \text{if } m_{s, t} =\infty.
\end{array}\right.
\end{equation}

Observe that 
\begin{equation}\label{eq:eps_st:2} 
\varepsilon_{s, t} =  
\displaystyle \left( \frac{q_s + q_t + 2\sqrt{q_s q_t}\cos(2\pi/m_{s, t})}{(q_s+1)(q_t+1)}\right)^{\frac 1 2}
\end{equation}	
if  $m_{s, t} \in \{ 3, 4, 6, 8\}$.

We define a symmetric bilinear form $\varphi_{\mathrm g}$ on the space $V_{\mathrm g} = \bigoplus_{s \in S} \KK \alpha_s$, called the \textbf{geometric form}, by setting 
$$\varphi_{\mathrm g}(\alpha_s, \alpha_t) = -\varepsilon_{s, t}$$
for all $s, t \in S$. By definition, it is a Hermitian form. It is canonically defined if  $(W, S)$ is \textbf{$2$-spherical} (i.e. $m_{s, t} < \infty$ for all $s, t \in S$). In general, it depends on the choice of the scalars $\mu_{s, t}$. 

Observe that if $q_s = 1$ for all $s$ and if $\mu_{s, t} =1$ for all $s, t$ with $m_{s, t} =\infty$, the form $\varphi_{\mathrm g}$  coincides with the Tits form. 

\medskip
For each subset $J \subset S$, we set $V_J = \langle \alpha_s \mid s \in J\rangle$ and denote  by $\varphi_J$ the restriction of  $\varphi_{\mathrm g}$ to $V_J \times V_J$. We say that $J$ is \textbf{spherical} if $W_J = \langle J \rangle$ is finite.
The following observation will be useful. 

\begin{lemma}\label{lem:pair-pos-def}
For each spherical pair $J = \{s, t\}\subseteq S$, the form $\varphi_J$ is positive definite. 
\end{lemma}
\begin{proof}
This is clear if $m_{s, t}=2$. If $m_{s, t} >2$, we observe that 
\begin{align*}
1-\varepsilon_{s, t}^2 & \geq \frac{(1+q_s)(1+q_t) - q_s - q_t - \sqrt{2 q_s q_t}}{(1+q_s)(1+q_t )}\\
&= \frac{1+q_s q_t -  \sqrt{2 q_s q_t}}{(1+q_s)(1+q_t )}\\
&> \frac{1+q_s q_t -  2\sqrt{ q_s q_t}}{(1+q_s)(1+q_t )}\\
&=\frac{(1-\sqrt{ q_s q_t})^2}{(1+q_s)(1+q_t )}
\end{align*}
It follows that the Gram matrix of   $\varphi_J$   has a positive determinant. The conclusion follows, since the trace of that matrix is equal to $2$.
\end{proof}

We shall now endow $(V_{\mathrm g}, \varphi_{\mathrm g})$ with the structure of a Hermitian $\KK_{\mathbf q}[W]$-module, as follows. For all $s \in S$, we set
$$E_s = \frac{1+e_s}{1+q_s}.$$
By the quadratic relations, we have $E_s^2 = E_s$. Moreover $E_s^* = E_s$. For all $s$, we define the linear operator
$$\widetilde{E_s} : V_{\mathrm g} \to V_{\mathrm g} : v \mapsto v-\varphi_{\mathrm g}(v, \alpha_s) \alpha_s.$$
By definition, this is the $\varphi_{\mathrm g}$-orthogonal projection on the hyperplane $\alpha_s^\perp$.

\begin{proposition}\label{prop:geometric-Hecke}
The assignments $E_s \mapsto  \widetilde{E_s}$ extend to a linear representation $\rho \colon \KK_{\mathbf q}[W] \to \End(V_{\mathrm g}, \varphi_{\mathrm g})$ of $*$-algebras. 
\end{proposition}

\begin{proof} 
We need to check that the operators $\rho(e_s) = (1+q_s)   \widetilde{E_s} - \mathrm{I}$ satisfy the defining  relations of the Hecke algebra. 

By definition, we have $ \widetilde{E_s}^2 =  \widetilde{E_s}$. It follows that $\rho(e_s)$ satisfies the quadratic relations. 

We shall now check the Artin relation for $s \neq t \in S$ with $m_{s, t} < \infty$. Set $J = \{s, t\}$. By Lemma~\ref{lem:pair-pos-def} the space $V$ decomposes as the orthogonal direct sum of $V_J$ and its $\varphi_{\mathrm g}$-orthogonal complement $V_J^\perp$. Observe that the restriction of $\widetilde{E_s}$ and $\widetilde{E_t}$ to $V_J^\perp$ is the identity operator. In particular the restriction to $\rho(e_s)$  and $\rho(e_t)$ to $V_J^\perp$ are equal to $q_s \mathrm{I}$ and $q_t \mathrm{I}$ respectively. Hence those restrictions satisfy the Artin relation. We now check that the restrictions of $\rho(e_s)$  and $\rho(e_t)$  to $V_J$ satisfy the Artin relations. 

Relative to the basis $(\alpha_s, \alpha_t)$, we have 
$$\rho(e_s) = \left(\begin{matrix}
-1 & (1+q_s)\varepsilon_{s, t}\\
0 & q_s
\end{matrix}\right) 
\qquad \text{and} \qquad
\rho(e_t) = \left(\begin{matrix}
q_t & 0\\
(1+q_t)\varepsilon_{s, t} & -1 
\end{matrix}\right).$$
 Hence
$$\rho(e_s)\rho(e_t) = \left(\begin{matrix}
-q_t + (1+q_s)(1+q_t)\varepsilon_{s, t}^2 &  - (1+q_s)\varepsilon_{s, t}\\
 q_s(1+q_t)\varepsilon_{s, t} & -q_s
\end{matrix}\right) $$
and
$$\rho(e_t)\rho(e_s) = \left(\begin{matrix}
-q_t & q_t(1+q_s)\varepsilon_{s, t}\\
 - (1+q_t)\varepsilon_{s, t}& -q_s + (1+q_s)(1+q_t)\varepsilon_{s, t}^2 
\end{matrix}\right).$$ 
Both matrices have trace $-q_s -q_t + (1+q_s)(1+q_t)\varepsilon_{s, t}^2 $ and determinant $q_s q_t$. The Artin relations follow from a direct computation using Cayley--Hamilton's theorem. This ensures that $\rho $ is indeed a linear representation of $\KK_{\mathbf q}[W]$. 

By definition we have $\varphi_{\mathrm g}(\widetilde{E_s}v, w)= \varphi_{\mathrm g}(v, \widetilde{E_s}w)$ for all $v, w \in  V_{\mathrm g}$. Since the Hecke algebra $\KK_{\mathbf q}[W]$ is generated by the $E_s$, we infer that $\rho$ is a representation of $*$-algebras.
\end{proof}

The pair $(V_\rg, \varphi_\rg)$ is thus a Hermitian 	$\KK_{\mathbf q}[W]$-module. We call it a \textbf{geometric Hecke module}. As mentioned above, it is canonically defined if $(W, S)$ is $2$-spherical, but not in general as it depends on the choice of the constants $\mu_{s, t}$ in (\ref{eq:eps_st}). 
	
The quadratic relations imply that, in any linear representation of the Hecke algebra, the only possible eigenvalues of the image of the generator $e_s$ are $q_s$ or $-1$. The \textbf{symmetric representation} is the one-dimensional representation that sends $e_s$ to $q_s$ for all $s$ (equivalently $T_s$ is mapped to $1$). The \textbf{alternating representation} is the one-dimensional representation that sends $e_s$ to $-1$ for all $s$.  

The following lemma is straightforward; we omit the proof.  

\begin{lemma}\label{lem:symmetric-submodule}
Let  $(V_{\mathrm g}, \varphi_{\mathrm g})$ be a geometric Hecke-module. A vector $v \in V_{\mathrm g}$ is fixed by $E_s$ for all $s \in S$ if and only if $v \in \mathrm{Rad}(\varphi_{\mathrm g})$. 

In particular, if $\varphi_{\mathrm g}$ is non-degenerate, then the symmetric representation is not contained in the geometric Hecke-module. \qed
\end{lemma}

We can now characterize the irreducibility of the geometric module. For that purpose, we define a graph $\Delta(\varphi_\rg)$ with vertex set $S$, where a pair $\{s, t\}$ of distinct vertices forms an edge if and only $\varepsilon_{s, t} \neq 0$. The following lemma is immediate from the definitions. We omit the proof. 

\begin{lemma}\label{lem:connected}
\begin{enumerate}[label=(\roman*)]
\item If $\Delta(\varphi_\rg)$ is connected, then $(W, S)$ is irreducible. 

\item If $(W, S)$ is irreducible, and if for all $s, t \in S$ with $m_{s, t} = \infty$, we have $\mu_{s, t} \neq 0$, then $\Delta(\varphi_\rg)$ is connected.  
\end{enumerate}
\end{lemma}

In particular, the graph $\Delta(\varphi_\rg)$ is connected if $(W, S)$ is irreducible and $2$-spherical. 

\begin{proposition}\label{prop:geom-irreducibility}
Suppose that $\Delta(\varphi_\rg)$ is connected. Then every proper submodule of $V_{\mathrm g}$ is contained in $\mathrm{Rad}( \varphi_{\mathrm g})$. 

In particular  the geometric Hecke module is irreducible if and only if   $\varphi_{\mathrm g}$ is non-degenerate. 
\end{proposition}
\begin{proof}
Let  $Y \leq V_{\mathrm g}$ be a submodule. 

\medskip
We first claim that if there exists $s \in S$ with  $\alpha_s \in Y$, then $Y = V_\rg$. 

Indeed, for any $t \in S$, the connectedness  of   $\Delta(\varphi_\rg)$ ensures that there is a sequence  $(s = s_0, s_1, s_2, \dots, s_n = t)$ of elements of $S$ such that $\varepsilon_{s_{i-1}, s_i} \neq 0$ for all $i$. It follows that if $\alpha_{s_{i-1}} \in Y$, then $\widetilde{E_{s_i}}(\alpha_{s_{i-1}}) \in Y$, hence $\alpha_{s_{i}} \in Y$. We deduce by induction that $\alpha_{s_{i}} \in Y$ for all  $i$, and the claim follows.

 \medskip
 Let now $s \in S$ and suppose that $Y$ contains some element  $y \not \in \alpha_s^\perp$. Then $\widetilde{E_{s}}(y) \in Y$, whence $\alpha_s \in Y$. By the claim, this implies that $Y = V_\rg$. 
 
It follows that every proper submodule of $V_\rg$ is contained in  $\bigcap_{s \in S}  \alpha_s^\perp = \mathrm{Rad}( \varphi_{\mathrm g})$, as required. 
\end{proof}

 \subsection{Signature of the geometric form}

The signature of the geometric form plays an important role. The following notation will be useful in studying it. We let $M_\qq(W)$ denote the Gram matrix of the geometric form $\varphi_\rg$ with respect to the standard basis $(\alpha_s)_{s \in S}$, and let $\lambda_\qq(W)$ denote the smallest eigenvalue of $M_\qq(W)$. Given two $|S|$-tuples $\qq= (q_s)_{s \in S}$ and $\mathbf r = (r_s)_{s \in S}$ of real numbers, we write $\qq \geq \mathbf r$ if $q_s \geq r_s$ for all $s \in S$. 

\begin{lemma}\label{lem:compare}
Let $(W, S)$ be a $2$-spherical Coxeter system of finite rank. Given  $S$-tuples $\qq= (q_s)_{s \in S}$ and $\mathbf r = (r_s)_{s \in S}$ of positive real numbers, if $q_s \geq r_s\geq 1$ for all $s \in S$, then  $\lambda_\qq(W) \geq \lambda_{\mathbf r}(W)$.  
\end{lemma}
\begin{proof}
All diagonal entries of the Gram matrix $M_\qq(W)$ are $1$. The off diagonal entry indexed by $(s,t)$ is $-\varepsilon_{s, t}$, where $\varepsilon_{s, t}$ is   defined by (\ref{eq:eps_st}). Observe that $\varepsilon_{s, t}$ is a non-increasing function of $q_s$ (resp. $q_t$). Hence the result follows from \cite[Proposition~2.9]{GriOpp}.
\end{proof}

\begin{proposition}\label{prop:spherical-affine}
Let $(W, S)$ be a Coxeter system of finite rank, and  $\qq= (q_s)_{s \in S}$ consist of real numbers.
\begin{enumerate}[label=(\roman*)]
\item If $(W, S)$ is spherical  and $q_s \geq 1$ for all $s$, then the geometric form is positive definite. 

\item If $(W, S)$ is  $2$-spherical affine  and $q_s > 1$ for all $s$, then the geometric form is positive definite. 
\end{enumerate}
	
\end{proposition}
\begin{proof}
 The Tits form of a spherical Coxeter system is positive definite, see \cite[Ch.~V, \S4.8]{Bour}. Hence (i) follows from Lemma~\ref{lem:compare}. 
 
 The Tits form of an   affine Coxeter system is positive semi-definite, with zero determinant,  see \cite[Ch.~V, \S4.9]{Bour}. Let $r = \min\{q_s |s \in S\}$, and let $\mathbf r$ be the constant $S$-tuple with $r_s =r$ for all $s$. By Lemma~\ref{lem:compare}, if the matrix $M_{\mathbf r}$  is positive definite, then $M_\qq(W)$ is positive definite as well. By hypothesis $r> 1$. The determinant of  $M_{\mathbf r}$ can be viewed as a rational function of $\sqrt r$, hence it has finitely many roots. Moreover, it vanishes at $r=1$, and that is non-negative for $r \geq 1$ by Lemma~\ref{lem:compare}. Hence it is strictly positive for $r>1$. The result follows. 
\end{proof}

\begin{corollary}\label{cor:lorentzian-dichotomy}
Let $(W, S)$ be a Coxeter system of finite rank $n = |S|$, and $\qq= (q_s)_{s \in S}$ consist of real numbers with $q_s > 1$ for all $s \in S$.

If $S$ contains a subset $J$ of size $n-1$ that is spherical or $2$-spherical affine, then  the Gram matrix $M_\qq(W)$ has at most one negative eigenvalue. In particular, the form $\varphi_\rg$ is positive-definite, or positive semi-definite and degenerate, or of Lorentzian signature and non-degenerate. 
\end{corollary}
\begin{proof}
Indeed the restriction $\varphi_J$ of the geometric form on the codimension one subspace $V_J$ is positive definite by Proposition~\ref{prop:spherical-affine}. Thus $V_\rg = V_J \perp V_J^\perp$ and $i_-(V_\rg) = i_-(V_J^\perp)$. The result follows. 
\end{proof}

\begin{remark}
In general, the signature of the geometric form can be of arbitrarily large   index even if $(W, S)$ is irreducible and $2$-spherical: indeed, if $S$ contains pairwise orthogonal subsets $J_1, \dots, J_r$ such that $\varphi_{J_i}$ has Lorentzian signature for all $i$, then the   index of $\varphi_\rg$ is at least $r$. Some examples will be discussed in \S\ref{sec:examples}. 
\end{remark}

\subsection{Kazhdan's property (T)}

Let us first recall Kassabov's criterion~\cite{Ka11}. 
Given a group $H$ generated by two subgroups $A, B \leq H$ 
 and a unitary representation $\pi$ of $H$ without non-zero $H$-invariant vectors, we define
$$\varepsilon_H(A, B; \pi) = \sup\big\{\frac{|\langle u, v \rangle |}{\| u \| \| v \|} \mid u \in \mathcal H_\pi^A \setminus \{0\}, v \in \mathcal H_\pi^B \setminus \{0\}\big\}.
$$
 If $\mathcal H_\pi^A = \{0\}$ or $\mathcal H_\pi^B = \{0\}$, we put $\varepsilon_H(A, B; \pi) =0$. If $\pi$ is an arbitrary unitary representation, we define  $\varepsilon_H(A, B; \pi)$ as $\varepsilon_H(A, B;\pi_0)$, where $\pi_0 \leq \pi$ is the sub-representation of $\pi$ defined on the orthogonal complement of the subspace of $H$-invariant vectors. 
The supremum of $\varepsilon_H(A, B; \pi) $ taken over all unitary representations $\pi$ of $H$ with $\mathcal H_\pi^H = \{0\}$, is denoted by 
$\varepsilon_H(A, B).$

Let us now consider a group $G$ generated by compact subgroups $A_1, \dots, A_n$, so that $A_{i, j} = \langle A_i \cup A_j \rangle$ is compact for all $i, j$. Set $\varepsilon_{i, j} = \varepsilon_{A_{i, j}}(A_i, A_j)$. Let $M$ be the symmetric $n\times n$-matrix defined by $M_{i, i} = 1$ and $M_{i, j}= -\varepsilon_{i, j}$  for all $i \neq j$. We call it the \textbf{cosine matrix} of the tuple $(G, (A_i)_{1\leq i \leq n})$. Kassabov's criterion \cite[\S2]{Ka11} ensures that if the cosine matrix $M$ is positive definite, then $G$ has property (T). 

Given a group  $G$ acting properly, chamber transitively on a locally finite building, the stabilizers of spherical residues are all compact. Moreover, the chamber-transitivity implies that $G$ is generated by the stabilizers of the panels containing any given chamber. The following result identifies the cosine matrix. 

\begin{proposition}\label{prop:cosine-matrix}
Let $X$ be a locally finite building of $2$-spherical type $(W, S)$ and $G$ be a locally compact group acting continuously, properly, chamber-transitively by  automorphisms on $X$. Let $c_0 \in X$ be a base chamber and for each $s\in S$, let $P_s \leq G$ be the stabilizer of the $s$-panel containing $c_0$. Let also $\qq=(q_s)_{s\in S}$ be the thickness vector of $X$. 

Then the cosine matrix of $(G, (P_s)_{s \in S})$ coincides with the Gram matrix $M_\qq(W)$ of the geometric form. 
\end{proposition}
\begin{proof}
Let $s \neq t \in S$. We denote by $P_{s, t}$ the subgroup generated by $P_s \cup P_t$. Since $G$ is chamber-transitive, the group $P_{s, t}$ coincides with the stabilizer of the $\{s, t\}$-residue of $X$ containing the base chamber $c_0$. 

Combining \cite[Proposition~7.1]{Garland} and \cite[Theorem~2.10(ii)]{CKKW}, we may now infer that   $\varepsilon_{P_{s, t}}(P_s, P_t) = \varepsilon_{s, t}$, where $\varepsilon_{s, t}$ is defined in (\ref{eq:eps_st}). The conclusion follows. 
\end{proof}

When the geometric form is positive definite, the induced form on the corresponding Hecke-induced module is positive semi-definite by Corollary~\ref{cor:main}. The corresponding unitary representation of $G$, obtained by dividing out the radical and completing, is then called the \textbf{geometric unitary representation} of $G$. 
 
 The following statement directly imples Corollary~\ref{corintro:Kazhdan} stated in the introduction. 
 
\begin{corollary}\label{cor:Kazhdan}
Let $X$, $G$ and $(P_s)_{s \in S}$ and $\qq$ be as in Proposition~\ref{prop:cosine-matrix}. Suppose that   $M_\qq(W)$ is positive definite. 

Then $G$ has Kazhdan's property (T). More precisely, the spectrum of the Laplacian operator defined by 
$$\Delta = \sum_{s\in S} (\one_B - (q_s+1)^{-1}\one_{P_s}) \in C^*(G)$$
is contained in $\{0\} \cup [\lambda_\qq(W), n]$, where $n = |S|$.  
Moreover, its spectral gap $\lambda_\qq(W)$ is realized as an eigenvalue of $\pi_\rg(\Delta)$, where $\pi_\rg$ is the  the geometric unitary representation of $G$. 

Furthermore,  if $(W, S)$ is irreducible and if  $G$ is Weyl-transitive, then the geometric unitary represention is  irreducible.
\end{corollary}
\begin{proof}
By Proposition~\ref{prop:cosine-matrix}, the cosine matrix $(G, (P_s)_{s \in S})$ coincides with the Gram matrix $M_\qq(W)$.  Since $G$ is generated by the parabolic subgroups $P_s$, which are compact, the positive definiteness of $M_\qq(W)$ ensures that $G$ has (T) by \cite[\S2]{Ka11}. 

For the other assertions of the corollary, we argue as in \cite[Theorem~6.1]{Ka11}. Let $\pi \colon G \to U(\mathcal{H})$ be a unitary representation. Observe first that the operator $\pi(\Delta)$ vanishes on $\mathcal H^G$ as well as on $(\mathcal H^B)^\perp$. 
It remains to describe the spectrum of $\pi(\Delta) $ on $\mathcal H^B \cap (\mathcal H^G)^\perp$. 

For each closed subspace $Z$ in $\mathcal H$, we denote by $d_Z(\cdot)$ the distance function to $Z$.   Given a vector $v\in \mathcal{H}^B$ that is orthogonal to $\mathcal{H}^G$, we  invoke \cite[Theorem~5.1]{Ka11}, which   ensures that  
$$ ||v||^2 = d_{\mathcal{H}^G}(v)^2 \leq D^T (M_{\textbf{q}}(W))^{-1} D , $$ 
where   $D$ is the column vector indexed by $S$ and defined by $D_s = d_{\mathcal{H}^{P_s}}(v)$. 

Let  $\lambda  = \lambda_\qq(W)> 0$. Using that $(q_s+1)^{-1}\pi(\one_{P_s})$ is the orthogonal projection on the subspace $\mathcal{H}^{P_s}$, we obtain 
\begin{align*}
D^T (M_{\textbf{q}}(W))^{-1} D & \leq \lambda^{-1} D^T \cdot D\\
 & = \lambda^{-1} \sum_{s\in S}d_{\mathcal{H}^{P_s}}(v)^2\\ 
& = \lambda^{-1} \sum_{s\in S}  \langle v - (q_s+1)^{-1}\pi(\one_{P_s}) v , v\rangle \\
& = \lambda^{-1} \sum_{s\in S} \langle (\pi(\one_B) - (q_s+1)^{-1} \one_{P_s}) v , v \rangle \\
&= \lambda^{-1} \langle \pi(\Delta)v,v \rangle, 
\end{align*}
where we have made use of the fact that $\pi(\textbf{1}_B)v = v$ because $v\in \mathcal{H}^B$. This proves that $\langle \pi(\Delta) v,v \rangle \geq \lambda ||v||^2$. Since $\Delta$ is self-adjoint, the spectrum of the operator $\pi(\Delta)$ consists of approximate eigenvalues (see \cite[Ch.~XI, \S1, Ex.~2]{Conway}). Therefore, we can conclude that the spectrum of $\pi(\Delta)$ is contained in $\{0\} \cup [\lambda, +\infty[$. 

Observe that the operator  $ \pi(\one_B - (1+q_s)^{-1} \one_{P_s} )$ is the orthogonal projection onto $\mathcal H^B \cap (\mathcal H^{P_s})^\perp$. In particular it has norm equal to $1$, so that $\|\pi(\Delta)\| \leq  |S|$. We infer that the spectrum of $\Delta$ in $C^*(G)$ is actually contained in $\{0\} \cup [\lambda,n]$, where $n = |S|$.

Finally, we observe  that the spectral gap is achieved in the geometric unitary representation. Indeed, let us denote this representation by $\pi_\rg \colon G \to U(\mathcal{K})$, and let $q \colon \Ind(V_0) \to \mathcal K$ be the canonical map. Let $c_0 \in X$ be the base chamber fixed by $B$, so that $q(\one_{c_0} \otimes V_0)$ is a $B$-invariant subspace of $\mathcal K$, which is isomorphic to $V_\rg$ since $q$ is injective on $\one_{c_0} \otimes V_0$ by the positive definiteness of $\varphi_\rg$. Since $G$ is chamber-transitive, the Hecke algebra $\CC_\qq[W]$ can be identified with a subalgebra of $C_c(B\backslash G/B)$, see  \S\ref{sec:Hecke-convol}.  Identifying  $q(\one_{c_0} \otimes V_0)$ with $V_\rg$, we see that $\pi_\rg(\Delta)$ acts on $V_{\rg}$ as $\sum_{s \in S} (1 - E_s) \in \CC_\qq[W]$. 
A direct computation  shows that the matrix representing $\sum_{s \in S} (1 - E_s)$ with respect to the standard basis $(\alpha_s)_{s \in S}$ is precisely  $M_{\textbf{q}}(W)$. Hence $\lambda$ is  an eigenvalue of $\pi_\rg(\Delta)$. 

The last statement on irreducibility   follows from Propositions~\ref{prop:Hecke-irred} and~\ref{prop:geom-irreducibility} if $(W, S)$ is irreducible and $G$ is Weyl-transitive. 
\end{proof}

\subsection{Isometric actions on infinite-dimensional symmetric spaces of finite rank}\label{sec:actions-sym-sp}

In general, the geometric form can be indefinite. The following result records that the geometric Hecke-module always satisfies the technical condition from Theorem~\ref{thm:main}. 

\begin{lemma}\label{lem:geom-form-tech}
Let $(W, S)$ be a Coxeter system of finite rank, and $\qq = (q_s)_{s \in S}$ be an $S$-tuple of  real numbers with $q_s \geq 1$ for all $s \in S$. Let $(V_\rg, \varphi_\rg)$ be the geometric $\KK_\qq[W]$-module. Suppose that $i_-(\varphi_\rg) = p$. For each spherical subset $J  \subset S$, there is a   $p$-dimensional subspace $V_J$ of $V_\rg$ that is pointwise fixed by $T_s$  for all $s \in J$, on which the restriction of $\varphi_\rg$ is negative definite. 
\end{lemma}
\begin{proof}
Consider the subspace $U_J$ spanned by $\{\alpha_s\mid s \in J\}$. The restriction $\varphi_J$ of the geometric form on  $U_J \times U_J$ is positive definite by Proposition~\ref{prop:spherical-affine}. Thus $V_\rg$ has an orthogonal sum decomposition as $U_J \perp U_J^\perp$, and $i_-(\varphi_J)=0$. In particular  $p = i_-(\varphi_\rg) = i_-(\varphi_\rg|_{U_J^\perp})$. By definition, the operator $T_s$ acts as the identity on $U_J^\perp$ for all $s \in J$. It follows that $U_J^\perp$ contains a subspace $V_J$ with the required properties. 
\end{proof}

Using Theorem~\ref{thm:main} this yields isometric actions of $\Aut(X)$ on  infinite-dimensional non-positively curved symmetric spaces of finite rank (see \S\ref{sec:sym-spaces} above for their definition).

\begin{theorem}\label{thm:action-sym-space}
Let $X$ be a locally finite semi-regular building of non-spherical type $(W, S)$ with thickness vector $\qq$, let $\KK$ be the field of real or complex numbers. We assume that every spherical residue of $X$ is chamber-transitive. 
Let $(V_\rg, \varphi_\rg)$ be the geometric $\KK_\qq[W]$-module defined above. Assume that $p = i_-(\varphi_\rg)>0$. 

\begin{enumerate}[label=(\roman*)]
\item The induced form $\Ind(\varphi_\rg)$ has negative index $p$, and its positive index, denoted by $q$, satisfies $q \geq 1$ (possibly $q = \infty$). 

\item The group $\Aut(X)$ has a continuous isometric action on the symmetric space $M_{p, q}(\KK)$ (see \S\ref{sec:sym-spaces}).
\end{enumerate}
Suppose in addition that $(W, S)$ is irreducible,  and moreover that  the graph $\Delta(\varphi_\rg)$ is connected (see Lemma~\ref{lem:connected}). Let $G$ be a locally compact group acting properly by automorphisms on $X$. 
\begin{enumerate}[resume*]
\item If $\Aut(X)$ is chamber-transitive and $G$ acts cocompactly on $X$,   then the $G$-action on $M_{p, q}(\KK)$ has no global fixed point. 

\item If $X$ is thick and $G$ acts Weyl-transitively on $X$, then $q = \infty$ and the $G$-action on $M_{p, q}(\KK) = X_p(\KK)$ is irreducible. 
\end{enumerate}
	
\end{theorem}
\begin{proof}
Let $V_0= V_\rg/\Rad(\varphi_\rg)$, and $\varphi_0$ be the form defined by $\varphi_\rg$ on the quotient. Hence $(V_0, \varphi_0)$ is a non-degenerate Hermitian $\KK_\qq[W]$-module. Let $p=i_-(V_0)$. Clearly $p = i_-(V_\rg)$. By Lemma~\ref{lem:geom-form-tech}, the technical condition in Theorem~\ref{thm:main} is satisfied. It follows that the negative index of $\Ind(\varphi_\rg)$ is equal to $p$. The fact that $i_+(\varphi_\rg) \geq 1$ is clear from the definition. Hence $q \geq 1$. This proves (i). 

The induced representation of $\Aut(X)$ on $V = \Ind(V_0)$ preserves the induced form $\varphi = \Ind(\varphi_0)$. By construction all vectors have an open stabilizer, hence the same holds for $V/\Rad(\varphi)$. Slightly abusing notation, we shall write $(\overline V, \overline \varphi)$  for the canonical completion of $V/\Rad(\varphi)$, equipped with the corresponding form. The associated representation on the  canonical completion $(\overline V, \overline \varphi)$ is also continuous   by Proposition~\ref{prop:continuity} (see \S\ref{sec:Hermitian}). Assertion (ii) follows. 

Assume now that $(W, S)$ is irreducible  and   that  $\Delta(\varphi_\rg)$ is connected. Suppose there is an $G$-fixed point in $M_{p, q}(\KK)$. Since $G$ is cocompact in $\Aut(X)$, it follows that $\Aut(X)$ has a bounded orbit, hence it also has a fixed because $M_{p, q}(\KK)$ is a complete CAT($0$) space. Thus we may assume without loss of generality that $G$ is chamber-transitive, since $\Aut(X)$ is chamber-transitive by hypothesis. 

By definition of $M_{p, q}(\KK)$,  there is a $p$-dimensional negative definite $G$-invariant subspace $W$ in $\overline V$. Choose a base chamber $c_0 \in X$. We have 
$$V_0 \cong \one_{c_0} \otimes V_0 \leq V.$$
Since $\varphi_0$ is non-degenerate, the projection $q \colon V \to \overline V$ is injective on $ \one_{c_0} \otimes V_0$. Moreover $G( \one_{c_0} \otimes V_0)$ spans $V$ by the hypothesis that $G$ is chamber-transitive. It follows that the image of $V_0 \cong \one_{c_0} \otimes V_0$  in $\overline V$ cannot be orthogonal to $W$. Thus the orthogonal projection of $W$ to $ q(\one_{c_0} \otimes V_0$) is non-zero. We denote by $W_0$ the corresponding subspace of $V_0$. Recall from \S\ref{sec:Hecke-convol} that $\KK_\qq[W]$ can be realized as a subalgebra of the convolution algebra $C_c^\KK(B\backslash G/B) \subseteq C_c^\KK(G)$, where $B = G_{c_0}$. Since $W$ is $G$-invariant, it is a $C_c^\KK(B\backslash G/B)$-submodule of $\overline V$.  It follows that  $W_0$ is a submodule of $V_0$. In view of Proposition~\ref{prop:geom-irreducibility}, we deduce that $W_0=V_0$. Hence 
$p \leq \dim(V_0) = \dim(W_0) \leq \dim(W) =p = i_-(V_0)$. Hence $V_0$ is negative definite. This is impossible by Proposition~\ref{prop:spherical-affine}. Assertion (iii) follows. 

Let us finally assume that $X$ is thick and that $G$ acts Weyl-transitively on $X$. In that case the natural embedding of $\KK_{\mathbf q}[W]$ in $C_c^\KK(B\backslash G/B)$ is surjective. The space $\overline V^B$ of $B$-invariant vectors contains $q( \one_{c_0} \otimes V_0)$. Moreover the Hecke-module $V_0$ is irreducible by Proposition~\ref{prop:geom-irreducibility}. Since $G \big(q(\one_c \otimes V_0)\big)$ spans a dense subspace of $\overline V$, we deduce from  Proposition~\ref{prop:Hecke-irred} that the $G$-representation to the isometry group of $M_{p, q}(\KK)$ is irreducible. If $q$ were finite, then  that group would be a finite-dimensional Lie group. Let $K \leq G$ be the kernel of that representation. By construction $K$ contains the kernel of the $G$-action on $X$, and the $G$-representation to the isometry group of $M_{p, q}(\KK)$ is defined on the whole group $\Aut(X)$, which is totally disconnected. Since Lie groups have no small subgroups, every continuous homomorphism from a totally disconnected locally compact group to a Lie group has an open kernel. We infer that $K$ is open in $G$. Using \cite[Proposition~4.4]{CapTho}, we deduce that $G/K$ is finite. Hence the $G$-orbits on the CAT($0$)-space $M_{p, q}(\KK)$ are finite. It follows that $G$ has a fixed point, which contradicts (iii). This confirms that $q =\infty$. 
\end{proof}

\begin{remark}
If the building $X$ is thin and $G = W$ (so that $G$ is Weyl-transitive), we have $\KK_\qq[W] \cong J = C_c^\KK(W)$, so that the induced module is nothing but the geometric module. It is thus finite-dimensional, thereby showing that the case $q< \infty$ can happen in Theorem~\ref{thm:action-sym-space}: the thickness hypothesis in the item (iv) is essential. The main case of interest is however when $X$ is thick and $\Aut(X)$ is Weyl-transitive, which is the case of the Kac--Moody buildings or the semi-regular right-angled buildings mentioned in the introduction. In all those cases $q =\infty$ by Theorem~\ref{thm:action-sym-space}(iv), even if the $G$-action is not Weyl-transitive.    
\end{remark}        

\begin{remark}
An irreducible representation afforded by Theorem~\ref{thm:action-sym-space}(iv) can fail to be \textit{geometrically irreducible}: indeed, an irreducible representation on $X_p(\CC)$ can preserve an isometrically embedded copy of $X_p(\RR)$. This is for example the case if the complex representation is the complexification of a real representation. We refer to  \cite[Proposition~7]{PySanchez} for an instance where every irreducible complex representation is of that kind.  
\end{remark}        

We derive the following, where the $G$-action is not assumed to be proper and $X$ is not assumed to be irreducible.

\begin{corollary}\label{cor:Lorentzian}
Let $X$ be a locally finite  building of non-spherical type $(W, S)$ with $\Aut(X)$ chamber-transitive, and  $\qq$ be its thickness vector. Let also  $G$ be a locally compact  group acting continuously, cocompactly on $X$

If $i_-(\varphi_\rg)=1$, then $G$ does not have Kazhdan's property (T). 
\end{corollary}
\begin{proof}
 By \cite[\S2]{DymJan}, there is no loss of generality in assuming that $(W, S)$ is $2$-spherical. Moreover $G$ acts continuously and cocompactly on every irreducible component of $X$. The hypotheses imply that there is an irreducible component $J \subset S$ such that the restriction of $\varphi_\rg$ to $V_J = \langle \alpha_s \mid s \in J \rangle$ has negative index~$1$. Let $X_J$ be the $J$-residue of a base chamber $c_0$. In view of Theorem~\ref{thm:action-sym-space}(iii), the closure of the natural image of $G$ in $\Aut(X_J)$ does not have (T) by Proposition~\ref{prop:FH}. Hence $G$ does not have (T) either. 
\end{proof}
	
\begin{proof}[Proof of Corollary~\ref{corintro:non-Kazhdan}]
In view of Proposition~\ref{prop:spherical-affine}, the hypotheses imply that $i_+(\varphi_\rg) \geq |S|-1$, so that $i_-(\varphi_\rg) \leq 1$. Since $\varphi_\rg$ is non-degenerate, we infer that either $\varphi_\rg$ is positive definite, or $i_-(\varphi_\rg) = 1$. The required conclusion follows from Corollaries~\ref{cor:Kazhdan} and~\ref{cor:Lorentzian}.
\end{proof}

\subsection{Examples}\label{sec:examples}

In this section, we illustrate the results above for three specific types. The case of Coxeter systems of rank $n$ containing a spherical or affine subset of rank~$n-1$ are covered by Corollary~\ref{corintro:non-Kazhdan}. We shall thus look at examples where the hypotheses of the latter do not apply. In each of the following examples, we consider a thick, locally finite, semi-regular building $X$ of type $(W, S)$ with thickness vector $\qq$, and  a locally compact group $G$ acting continuously, properly, chamber-transitively by  automorphisms on $X$.

\begin{example}\label{ex:1}
The first example is that of a $2$-spherical Coxeter system $(W, S)$ of rank~$7$ with the following Coxeter diagram:
\begin{center}
\begin{tikzpicture}
\filldraw[black] (0,0) circle (2pt) ;
\filldraw[black] (1,0) circle (2pt) ;
\filldraw[black] (-1,0) circle (2pt) ;
\draw[black, thick] (0,0) -- (1,0);
\draw[black, thick] (-1,0) -- (0,0);
\filldraw[black] (-2,-.75) circle (2pt) ;
\filldraw[black] (-2,.75) circle (2pt) ;
\filldraw[black] (2,-.75) circle (2pt) ;
\filldraw[black] (2,.75) circle (2pt) ;
\draw[black, thick] (-2,-.75) -- (-1,0);
\draw[black, thick] (-2,.75) -- (-1,0);
\draw[black, thick] (2,-.75) -- (1,0);
\draw[black, thick] (2,.75) -- (1,0);
\draw[black, thick] (-2,-.75) -- (-2,.75);
\draw[black, thick] (2,-.75) -- (2,.75);
\end{tikzpicture}
\end{center}
For all generators $s, t \in S$ such that $m_{s, t}$ is an odd integer, the elements $s$ and $t$ are conjugate. We infer that all Coxeter generators are conjugate in the example at hand, so that the thickness vector $\qq=(q_s)_{s \in S}$ is constant. Let $q= q_s $ be the common value. 

Observe that $S$ contains a spherical subset of size $5$. Proposition~\ref{prop:spherical-affine} ensures that $i_+(\varphi_\rg) \geq 5$, so that $i_-(\varphi_\rg) \leq 2$. Set 
$$z = \sqrt q \qquad \text{and} \qquad
f(z) = (1+q)^7 \det(M_\qq(W)).$$ 
We   obtain 
\begin{align*}
f(z) 
& =z^{14} - z^{12} - 4z^{11} - 2z^{10} + 4z^8 + 4z^7 + 4z^6 - 2z^4 - 4z^3 - z^2 + 1\\
&=(z^2+z+1)^3 (z-  1 )^2 (z^6 - z^5 - z^4 - z^2 - z + 1).
\end{align*}
One checks that $f(z)<0$ for $z \in   \{\sqrt 2, \sqrt 3\}$ and $f(z)>0$ for all $z \geq 2$. It follows that the signature of  index of $M_\qq(W)$ is $(1, 6)$  for $q=2, 3$, and $(0, 7)$ for $q \geq 4$. In view of Corollaries~\ref{cor:Kazhdan} and~\ref{cor:Lorentzian}, we infer  that $G$ has (T) for $q \geq 4$, whereas $G$ does not have (T) for $q =2, 3$.
\end{example}

\begin{example}\label{ex:2}
The second example is that of a $2$-spherical  Coxeter system $(W, S)$ of rank~$7$ with the following Coxeter diagram:
\begin{center}
\begin{tikzpicture}
\filldraw[black] (0,0) circle (2pt) ;
\filldraw[black] (1,0) circle (2pt) ;
\filldraw[black] (-1,0) circle (2pt) ;
\draw[black, thick] (0,0) -- (1,0);
\draw[black, thick] (-1,0) -- (0,0);
\filldraw[black] (-2,-.75) circle (2pt) ;
\filldraw[black] (-2,.75) circle (2pt) ;
\filldraw[black] (2,-.75) circle (2pt) ;
\filldraw[black] (2,.75) circle (2pt) ;
\draw[black, thick] (-2,-.75) -- (-1,0);
\draw[black, thick] (-2,.75) -- (-1,0);
\draw[black, thick] (2,-.75) -- (1,0);
\draw[black, thick] (2,.75) -- (1,0);
\draw[black, thick] (-2,-.75) -- (-2,.75);
\draw[black, thick] (2,-.75) -- (2,.75);
\node at (2.2,0) {$4$};
\node at (-2.2,0) {$4$};
\end{tikzpicture}
\end{center}
As in Example~\ref{ex:1}, the thickness vector is constant; we let $q = q_s$ be the common value and $z = \sqrt q$. The set $S$ contains again a spherical subset of size $5$, so that  $i_-(\varphi_\rg) \leq 2$. Setting again $ f(z) = (1+q)^7 \det(M_\qq(W)) $, we   obtain 
\begin{align*}
f(z) 
& = z^{14} - 3z^{12} - 4\sqrt 2 z^{11} - z^{10} + 4 \sqrt 2 z^9 + 11 z^8 + 8 \sqrt 2 z^7 \\
& \qquad  + 11 z^6 + 4 \sqrt 2 z^5 - z^4 - 4 \sqrt 2 z^3 - 3 z^2 + 1\\
&=(z^2+ \sqrt 2 z+1)^2 (z^4 - \sqrt 2 z^3 - \sqrt 2 z + 1 ) (z^6 - \sqrt 2 z^5 - z^4 - z^2 - \sqrt 2 z + 1).
\end{align*}
One checks that $f(z)> 0$ for $z =\sqrt 2$ and all $z \geq \sqrt 5$, whereas $f(z)<0$ for $z \in \{\sqrt 3, 2\}$. Since the determinant of the $3\times 3$-submatrix associated with the compact hyperbolic Coxeter group of type $(3, 3, 4)$ is negative for $z=\sqrt 2$, it follows that $i_-(\varphi_\rg) \geq 1$ for $q=2$. We infer  that the signature of  index of $M_\qq(W)$ is $(2, 5)$  for $q=2$, $(1, 6)$ for $q=3, 4$, and $(0, 7)$ for $q \geq 5$. In view of Corollaries~\ref{cor:Kazhdan} and~\ref{cor:Lorentzian}, we infer  that $G$ has (T) for $q \geq 5$, whereas $G$ does not have (T) for $q =3, 4$. For $q =2$, Theorem~\ref{thm:action-sym-space} yields a representation of $G$ on an infinite-dimensional symmetric space of rank~$2$, which is irreducible if $G$ is Weyl-transitive. This action does not yield any conclusion regarding property (T). However, one checks from the presentation of the Hecke algebra $\RR_\qq[W]$ that it admits, as a quotient, the algebra $\RR_\qq[W']$, where $(W', S)$ is the Coxeter system of type 
\begin{center}
\begin{tikzpicture}
\filldraw[black] (0,0) circle (2pt) ;
\filldraw[black] (1,0) circle (2pt) ;
\filldraw[black] (-1,0) circle (2pt) ;
\draw[black, thick] (0,0) -- (1,0);
\draw[black, thick] (-1,0) -- (0,0);
\filldraw[black] (-2,-.75) circle (2pt) ;
\filldraw[black] (-2,.75) circle (2pt) ;
\filldraw[black] (2,-.75) circle (2pt) ;
\filldraw[black] (2,.75) circle (2pt) ;
\draw[black, thick] (-2,-.75) -- (-1,0);
\draw[black, thick] (-2,.75) -- (-1,0);
\draw[black, thick] (2,-.75) -- (1,0);
\draw[black, thick] (2,.75) -- (1,0);
\draw[black, thick] (2,-.75) -- (2,.75);
\node at (2.2,0) {$4$};
\end{tikzpicture}
\end{center}
The latter has a spherical subset of size~$6$, and one computes that the determinant of $M_\qq(W')$ is negative for $q=2, 3$. Argueing as in  Theorem~\ref{thm:main}(iii), we obtain a fixed-point-free action of $G$ on the infinite-dimensional real hyperbolic space, thereby showing that $G$ does not have (T) for $q=2$ either, by Proposition~\ref{prop:FH}.
\end{example}

\begin{example}\label{ex:3}
The third example is that of a right-angled  Coxeter system $(W, S)$ of rank~$4$ with the following Coxeter diagram:
\begin{center}
\begin{tikzpicture}
\filldraw[black] (0,0) circle (2pt) ;
\filldraw[black] (1,0) circle (2pt) ;
\filldraw[black] (2,0) circle (2pt) ;
\filldraw[black] (-1,0) circle (2pt) ;
\draw[black, thick] (0,0) -- (1,0);
\draw[black, thick] (2,0) -- (1,0);
\draw[black, thick] (-1,0) -- (0,0);
\node at (0.5,.3) {$\infty$};
\node at (-0.5,.3) {$\infty$};
\node at (1.5, .3) {$\infty$};
\node at (-1,-.3) {$1$};
\node at (0,-.3) {$2$};
\node at (1, -.3) {$3$};
\node at (2, -.3) {$4$};
\end{tikzpicture}
\end{center}
That group $W$ is Gromov hyperbolic (see \cite[\S12.6]{Davis}), so that $X$ and $G$ are hyperbolic as well. 
Set $a = \mu_{1, 2}$, $b = \mu_{2, 3}$ and $c = \mu_{3, 4}$, so that $$M_\qq(W) = \left(\begin{array}{cccc}
1 & -a & 0 & 0\\
- \bar a & 1 & -b & 0\\
0 &  - \bar b & 1 & -c\\
0 & 0 & - \bar c & 1
\end{array}\right)$$
and 
$$\det(M_\qq(W) ) = |a|^2 |c|^2 - |a|^2 - |b|^2 - |c|^2 + 1= (|a|^2-1)(|c|^2-1) - |b|^2.$$
Since $S$ contains   a spherical subset of size $2$, we have  $i_-(\varphi_\rg) \leq 2$. Moreover, if $|a| > 1$, or $|b|> 1$, or $|c|>1$, then $M_\qq(W)$ contains a $2\times2$-submatrix of Lorentzian signature, so that $i_-(\varphi_\rg) \geq 1$. In particular, for all $a$ and $c$ with $|a|, |c| > 1$, the matrix $M_\qq(W)$ has signature $(2, 2)$ for all $b$ with $0< |b| < \sqrt{(|a|^2-1)(|c|^2-1)}$. Thus Theorem~\ref{thm:action-sym-space} yields a multi-parameter set of representations of $G$ on  the real and complex infinite-dimensional symmetric space of rank~$2$.
\end{example}

\subsection{Metrically proper actions}\label{sec:proper}

In this action, we discuss whether the actions on infinite-dimensional symmetric spaces afforded by Theorem~\ref{thm:action-sym-space} are metrically proper. In case of the real hyperbolic space, that condition implies the Haagerup property, see Proposition~\ref{prop:FH}.

\begin{lemma}\label{lem:char-proper}
Let $X$ be a locally finite semi-regular building of non-spherical type $(W, S)$ with thickness vector $\qq$.  Suppose that $\Aut(X)$ is chamber-transitive,  and that all coefficients of $M_\qq(W)$ are real.  

Suppose that $i_-(\varphi_\rg)=1$. Then the $\Aut(X)$-action on the hyperbolic space $\mathbb H  = X_1(\RR)$ afforded by Theorem~\ref{thm:action-sym-space} is metrically proper if and only if for some negative vector $v_0 \in V_\rg$, the map $W \to \RR : w \mapsto \varphi_\rg(T_w v_0, v_0)$ is  proper. 

If these conditions hold, then every  closed subgroup of $\Aut(X)$ has the Haagerup property. 
\end{lemma}
\begin{proof}
We adopt the same notation as in the proof of  Theorem~\ref{thm:action-sym-space}, with   $\KK  =\RR$, and let $\sigma$ denote the resulting representation of $G = \Aut(X)$ on $\overline V$, preserving the  form $\overline\varphi$ that has negative index one.

Fix a base chamber $c_0\in X$, let $\xi$ be the image of $\one_{c_0}\otimes v_0$ in $\overline V$, and let $x_0=[\xi]$
be the corresponding point of real hyperbolic space. We have $\overline\varphi(\xi,\xi)=-1$. With the standard normalization of the hyperbolic distance $d_{\mathbb H}$, the definition of the induced form gives
$$
\cosh d_{\mathbb H}(x_0,gx_0)
=|\overline\varphi(\sigma(g)\xi,\xi)|=|\varphi_\rg(T_{\delta(c_0,gc_0)}v,v)|
$$
for all $g \in G$. 
Since $G=\Aut(X)$ acts properly on $X$, the map $G \to W : g \mapsto \delta(c_0,gc_0)$ is proper. The required equivalence    on metric properness follows. Those conditions imply that $\Aut(X)$, hence each of its   closed subgroups, has the Haagerup property by Proposition~\ref{prop:FH}(2).
\end{proof}

The following observation provides an obstruction to the properness of the map in Lemma~\ref{lem:char-proper}. 

\begin{lemma}\label{lem:non-proper}
Retain the notation and assumptions of Lemma~\ref{lem:char-proper}, and suppose that $i_-(\varphi_\rg)=1$.  If for some non-spherical subset $J \subseteq S$, the form $\varphi_J$ (defined as in \S\ref{sec:geom-Hecke-module}) is positive definite,   the $\Aut(X)$-action on the hyperbolic space $\mathbb H  = X_1(\RR)$ afforded by Theorem~\ref{thm:action-sym-space} is not metrically proper. 
\end{lemma}
\begin{proof}
The hypothesis implies that $\varphi_J$ is non-degenerate, so that $V_\rg$ has an orthogonal decomposition as $V_J \perp V_J^\perp$. Since $\varphi_\rg$ has negative index~$1$, there exists a negative vector $v_0 \in V_J^\perp$. By definition of the geometric module, we have $E_s(v_0) = v_0$ for all $s \in J$, hence $T_w(v_0) = v_0$ for all $w \in W_J$. Thus the restriction of the map $w \mapsto \varphi_\rg(T_w v_0, v_0)$ to the infinite subgroup $W_J$ is constant. Hence it is not proper, and the conclusion follows from Lemma~\ref{lem:char-proper}.
\end{proof}

The following result may be viewed as a converse to Lemma~\ref{lem:non-proper}. 

\begin{proposition}\label{prop:proper}
Let $X$ be a thick, locally finite semi-regular building of non-spherical type $(W, S)$ with thickness vector $\qq$.  Suppose that $\Aut(X)$ is chamber-transitive,  and that all coefficients of $M_\qq(W)$ are real.
Suppose moreover that for each  non-spherical subset $J \subseteq S$, the form $\varphi_J$ is non-degenerate with negative index~$1$. Then the $\Aut(X)$-action on the hyperbolic space $\mathbb H  = X_1(\RR)$ afforded by Theorem~\ref{thm:action-sym-space} is metrically proper. 
\end{proposition}
\begin{proof}
The formula defining the geometric Hecke module gives
$$
T_sx=x-(1+q_s^{-1})\varphi_\rg(x,\alpha_s)\alpha_s
$$
for all $x \in V_\rg$ and $s \in S$. 
Since $\varphi_\rg(\alpha_s,\alpha_s)=1$, a direct computation
yields
\begin{equation}\label{eq:Haagerup-generator}
\varphi_\rg(T_sx,T_sx)
=
\varphi_\rg(x,x)
-(1-q_s^{-2})|\varphi_\rg(x,\alpha_s)|^2.
\end{equation}

Given a reduced expression $w = s_n \dots s_1$, we consider an operator
$A=T_{s_n}\cdots T_{s_1}$, and  we define  $A_0=\mathrm{Id}$ and $A_j=T_{s_j}\cdots T_{s_1}$ for
$1\leq j\leq n$. By using  
\eqref{eq:Haagerup-generator} repeatedly, we obtain
\begin{equation}\label{eq:Haagerup-block}
\varphi_\rg(x,x)-\varphi_\rg(Ax,Ax)
=\sum_{j=1}^{n}(1-q_{s_j}^{-2})
  |\varphi_\rg(A_{j-1}x,\alpha_{s_j})|^2.
\end{equation}
Since $X$ is thick, we have $q_s \geq 2$ for all $s$. Therefore
\begin{equation}\label{eq:basic}
\varphi_\rg(x,x) \geq \varphi_\rg(T_w x,T_w x)  
\end{equation}
for all $x \in V_\rg$ and $w \in W$.

The \textbf{support} of an element $w \in W$, denoted by $\supp(w)$, is defined as the smallest subset $I \subseteq S$ such that $w \in W_I$. Since $W_I \cap W_J = W_{I \cap J}$, the support is well-defined. Let $L$ be the maximal length of an element of $W$ with spherical support (notice that $L < \infty$ since $S$ is finite), and  set $N = L+1$. 

\medskip
We claim that there exist a constant $\eta > 1$ such that for each non-spherical subset $J \subseteq S$, each $w \in W$ with $\ell(w)=N$ and $\supp(w) = J$, and each $x \in V_J$, we have 
$$-\varphi_\rg(T_w x,T_wx)
\geq \eta\bigl(-\varphi_\rg(x,x)\bigr).$$

Indeed, consider a reduced expression  $w =s_N \dots s_1$. Set $A = T_w$ and $Q_A(x) = \varphi_\rg(x,x)-\varphi_\rg(Ax,Ax)$. Given $x \in V_J$ with $Q_A(x) = 0$. Then every summand  in (\ref{eq:Haagerup-generator}) vanishes. The formula for $T_{s_j}$ then implies successively that
$A_jx=x$ for all $j$. Consequently,
$\varphi_\rg(x,\alpha_{s_j})=0$ for every $j$.
By definition all elements of $J$ occur, and by assumption we have $x \in V_J$. The hypothesis that  $\varphi_J$ is non-degenerate
therefore implies that $x=0$. This shows that $Q_A$ is a positive-definite quadratic form on $V_J$. 
Since the number of elements $w$ of length $N$ in $W$ is finite, and since   the cone of positive definite forms on the finite-dimensional
space $V_J$ is open in the space of symmetric bilinear forms on $V_J$, there exists $\varepsilon_J>0$ such that
$
Q_{T_w}(x)+\varepsilon \varphi_\rg(x,x)\geq0
$
for every $w$ with $\ell(w)=N$ and $\supp(w) = J$,  every $x\in V_J$, and every $\varepsilon \in [0, \varepsilon_J]$.
Define $\varepsilon$ as the minimum of $ \varepsilon_J$ taken over all non-spherical subsets $J$. Then the constant  $\eta=1+\varepsilon>1$ satisfies
$$-\varphi_\rg(T_w x,T_w x)
=-\varphi_\rg(x,x)+Q_{T_w}(x)
\geq \eta\bigl(-\varphi_\rg(x,x)\bigr)
$$
for all $w$  with $\ell(w)=N$ and $\supp(w) = J$, and all $x\in V_J$.  This proves the claim. 

\medskip
Let now $w \in W$ with $\ell(w)=N$ and $\supp(w) = J$, and let $x \in V_\rg$. We do not assume that $x \in V_J$. By hypothesis $\varphi_J$ is non-degenerate, so we have an orthogonal decomposition $V_\rg = V_J \perp V_J^\perp$. Let $x = x_J + x'$ be the corresponding decomposition of $x$. Since $x' \in V_J^\perp$, we have $T_s(x') = x'$ for all $s \in J$. In particular $T_w(x') = x'$. In view of the claim, we obtain
\begin{align*}
-\varphi_\rg(T_w x,T_w x) 
&= -\varphi_\rg(T_w x_J + x',T_w x_J+x')\\
&= -\varphi_\rg(T_w x_J,T_w x_J) -\varphi_\rg( x',x') \\
&\geq \eta\bigl(-\varphi_\rg(x_J,x_J)\bigr)  -\varphi_\rg( x',x') \\
&=  \eta\bigl(-\varphi_\rg(x,x) + \varphi_\rg(x', x')\bigr)  -\varphi_\rg( x',x') \\
&=  \eta\bigl(-\varphi_\rg(x,x) \bigr) + (\eta-1) \varphi_\rg( x',x').
\end{align*}
By hypothesis, the forms $\varphi_J$ and $\varphi_\rg$ both have negative index one. It follows that the restriction of $\varphi_\rg$ on $V_J^\perp$ is positive-definite. Moreover $\eta>1$. Hence the term $ (\eta-1) \varphi_\rg( x',x')$ is non-negative, and we conclude that the inequality 
$$
-\varphi_\rg(T_w x,T_w x) \geq  \eta\bigl(-\varphi_\rg(x,x) \bigr) 
$$
holds for all $x \in V_\rg$.

Choose $v_0 \in V_\rg$ with $\varphi_\rg(v_0, v_0) = -1$. Given $w \in W$, we may decompose it   as a product of $\lfloor\ell(w)/N\rfloor$ elements of length $N$, times an element $w'$ of length~$<N$.  Invoking the previous inequality repeatedly, together with (\ref{eq:basic}), 
we obtain
\begin{align*}
-\varphi_\rg(T_wv_0,T_wv_0)
& \geq
\eta^{\lfloor\ell(w)/N\rfloor}\bigl(-\varphi_\rg(T_{w'}v_0, T_{w'}v_0)\bigr) \\
&\geq \eta^{\lfloor\ell(w)/N\rfloor}\bigl(-\varphi_\rg(v_0, v_0)\bigr) \\
&=  \eta^{\lfloor\ell(w)/N\rfloor}.
\end{align*}

For any vector $x \in V_\rg$, let $x = \lambda v_0 + x'$ be the decomposition of $x$ according to $V_\rg = \langle v_0 \rangle \perp v_0^\perp$. We have $\lambda = -\varphi_\rg(x, v_0)$ and $\varphi_\rg(x, x) = -|\lambda|^2 + \varphi_\rg(x', x') \geq - |\lambda|^2$, so that 
$$\varphi_\rg(x, x)  \geq - |\varphi_\rg(x, v_0)|^2.$$
For $x = T_wv_0$, we deduce  that $|\varphi_\rg(T_w v_0, v_0)| \geq  \eta^{\frac 1 2 \lfloor\ell(w)/N\rfloor}$. This confirms that the map $w \mapsto \varphi_\rg(T_wv_0, v_0)$ is proper. The conclusion follows in view of Lemma~\ref{lem:char-proper}.
\end{proof}

\begin{proof}[Proof of Corollary~\ref{corintro:Haagerup}]
By hypothesis, every proper subset $J \subsetneq S$ is spherical, and $\varphi_\rg$ has negative index~$1$. Hence $\varphi_\rg$  is  non-degenerate by Corollary~\ref{cor:lorentzian-dichotomy}. We infer that the hypotheses of Proposition~\ref{prop:proper} are fulfilled. The conclusions follow, in view of Lemma~\ref{lem:char-proper}.
\end{proof}		

Let us  mention that Proposition~\ref{prop:proper} also applies to groups that are not of compact hyperbolic type. 

\begin{example}
Let $X$ be a thick, locally finite, semi-regular building with $\Aut(X)$ chamber-transitive, and whose  type is the Coxeter system $(W', S)$ appearing in Example~\ref{ex:2}. If $q_s =2$ for all $s$, then the $\Aut(X)$-action on the real hyperbolic space afforded by Theorem~\ref{thm:action-sym-space} is metrically proper (hence $\Aut(X)$ has the Haagerup property). Indeed,  one checks that the hypothesis of Proposition~\ref{prop:proper} is satisfied. This is no longer the case if $q_s =3$ for all $s$: indeed, if $J \subset S$ is the subset of compact hyperbolic type, the form $\varphi_J$ is positive definite in that case, so that the hyperbolic action is not proper by Lemma~\ref{lem:non-proper}.
\end{example}

\bibliographystyle{abbrv} 
\bibliography{Hecke}

@book{CCJJV,
  author    = {Cherix, Pierre-Alain and Cowling, Michael and
               Jolissaint, Paul and Julg, Pierre and Valette, Alain},
  title     = {Groups with the {Haagerup} Property:
               {Gromov}'s a-{T}-menability},
  series    = {Progress in Mathematics},
  volume    = {197},
  publisher = {Birkh{\"a}user Verlag},
  address   = {Basel},
  year      = {2001}
}

@article {Ka11,
    AUTHOR = {Kassabov, Martin},
     TITLE = {Subspace arrangements and property {T}},
   JOURNAL = {Groups Geom. Dyn.},
  FJOURNAL = {Groups, Geometry, and Dynamics},
    VOLUME = {5},
      YEAR = {2011},
    NUMBER = {2},
     PAGES = {445--477},
      ISSN = {1661-7207,1661-7215},
   MRCLASS = {20F65 (20F55 22D10)},
  MRNUMBER = {2782180},
MRREVIEWER = {Pierre-Emmanuel\ Caprace},
       DOI = {10.4171/GGD/134},
       URL = {https://doi.org/10.4171/GGD/134},
}

@book {AbramenkoBrown,
    AUTHOR = {Abramenko, Peter and Brown, Kenneth S.},
     TITLE = {Buildings},
    SERIES = {Graduate Texts in Mathematics},
    VOLUME = {248},
      NOTE = {Theory and applications},
 PUBLISHER = {Springer, New York},
      YEAR = {2008},
     PAGES = {xxii+747},
      ISBN = {978-0-387-78834-0},
   MRCLASS = {20E42 (20F55 20J06 51E24 51F15)},
  MRNUMBER = {2439729},
MRREVIEWER = {Ralf\ Koehl},
       DOI = {10.1007/978-0-387-78835-7},
       URL = {https://doi.org/10.1007/978-0-387-78835-7},
}

@book {BHV,
    AUTHOR = {Bekka, Bachir and de la Harpe, Pierre and Valette, Alain},
     TITLE = {Kazhdan's property ({T})},
    SERIES = {New Mathematical Monographs},
    VOLUME = {11},
 PUBLISHER = {Cambridge University Press, Cambridge},
      YEAR = {2008},
     PAGES = {xiv+472},
      ISBN = {978-0-521-88720-5},
   MRCLASS = {22-02 (22E40 28D15 37A15 43A07 43A35)},
  MRNUMBER = {2415834},
MRREVIEWER = {Markus\ Neuhauser},
       DOI = {10.1017/CBO9780511542749},
       URL = {https://doi.org/10.1017/CBO9780511542749},
}

@article {Oppenheim,
    AUTHOR = {Oppenheim, Izhar},
     TITLE = {Property ({T}) for groups acting on affine buildings},
   JOURNAL = {Bull. Lond. Math. Soc.},
  FJOURNAL = {Bulletin of the London Mathematical Society},
    VOLUME = {57},
      YEAR = {2025},
    NUMBER = {10},
     PAGES = {3151--3162},
      ISSN = {0024-6093,1469-2120},
   MRCLASS = {20F65 (22D55)},
  MRNUMBER = {4971894},
       DOI = {10.1112/blms.70148},
       URL = {https://doi.org/10.1112/blms.70148},
}

@book {Davis,
    AUTHOR = {Davis, Michael W.},
     TITLE = {The geometry and topology of {C}oxeter groups},
    SERIES = {London Mathematical Society Monographs Series},
    VOLUME = {32},
 PUBLISHER = {Princeton University Press, Princeton, NJ},
      YEAR = {2008},
     PAGES = {xvi+584},
      ISBN = {978-0-691-13138-2; 0-691-13138-4},
   MRCLASS = {20F55 (05B45 05C25 51-02 57M07)},
  MRNUMBER = {2360474},
MRREVIEWER = {Ralf\ Koehl},
}

@book {Bour,
    AUTHOR = {Bourbaki, N.},
     TITLE = {\'El\'ements de math\'ematique. {F}asc. {XXXIV}. {G}roupes et
              alg\`ebres de {L}ie. {C}hapitre {IV}: {G}roupes de {C}oxeter
              et syst\`emes de {T}its. {C}hapitre {V}: {G}roupes engendr\'es
              par des r\'eflexions. {C}hapitre {VI}: syst\`emes de racines},
    SERIES = {Actualit\'es Scientifiques et Industrielles [Current
              Scientific and Industrial Topics]},
    VOLUME = {No. 1337},
 PUBLISHER = {Hermann, Paris},
      YEAR = {1968},
     PAGES = {288 pp. (loose errata)},
   MRCLASS = {22.50 (17.00)},
  MRNUMBER = {240238},
MRREVIEWER = {G.\ B.\ Seligman},
}

@article {DymJan,
    AUTHOR = {Dymara, Jan and Januszkiewicz, Tadeusz},
     TITLE = {Cohomology of buildings and their automorphism groups},
   JOURNAL = {Invent. Math.},
  FJOURNAL = {Inventiones Mathematicae},
    VOLUME = {150},
      YEAR = {2002},
    NUMBER = {3},
     PAGES = {579--627},
      ISSN = {0020-9910,1432-1297},
   MRCLASS = {20E42 (20F55 20J06 22E50)},
  MRNUMBER = {1946553},
MRREVIEWER = {Alain\ Valette},
       DOI = {10.1007/s00222-002-0242-y},
       URL = {https://doi.org/10.1007/s00222-002-0242-y},
}

@article {Bischof,
    AUTHOR = {Bischof, Sebastian},
     TITLE = {On growth functions of {C}oxeter groups},
   JOURNAL = {Proc. Edinb. Math. Soc. (2)},
  FJOURNAL = {Proceedings of the Edinburgh Mathematical Society. Series II},
    VOLUME = {68},
      YEAR = {2025},
    NUMBER = {3},
     PAGES = {979--993},
      ISSN = {0013-0915,1464-3839},
   MRCLASS = {20F55 (20F69 51F15)},
  MRNUMBER = {4943977},
MRREVIEWER = {Anthony\ Genevois},
       DOI = {10.1017/S0013091525000094},
       URL = {https://doi.org/10.1017/S0013091525000094},
}

@phdthesis{BischofThesis,
  author = {Bischof, Sebastian},
  title  = {Construction of {RGD}-systems of type $(4,4,4)$ over
            ${\mathbb F_2}$},
  school = {Justus-Liebig-Universit{\"a}t Gie{\ss}en},
  year   = {2023},
  url    = {https://jlupub.ub.uni-giessen.de/items/4e38869c-b9ea-4c74-a9fb-beb2056d2e6f}
}

@book {Tits74,
    AUTHOR = {Tits, Jacques},
     TITLE = {Buildings of spherical type and finite {BN}-pairs},
    SERIES = {Lecture Notes in Mathematics},
    VOLUME = {Vol. 386},
 PUBLISHER = {Springer-Verlag, Berlin-New York},
      YEAR = {1974},
     PAGES = {x+299},
   MRCLASS = {20G15 (50A20)},
  MRNUMBER = {470099},
MRREVIEWER = {Bruce\ Cooperstein},
}

@article {CKKW,
    AUTHOR = {Caprace, Pierre-Emmanuel and Conder, Marston and Kaluba, Marek
              and Witzel, Stefan},
     TITLE = {Hyperbolic generalized triangle groups, property ({T}) and
              finite simple quotients},
   JOURNAL = {J. Lond. Math. Soc. (2)},
  FJOURNAL = {Journal of the London Mathematical Society. Second Series},
    VOLUME = {106},
      YEAR = {2022},
    NUMBER = {4},
     PAGES = {3577--3637},
      ISSN = {0024-6107,1469-7750},
   MRCLASS = {20F67 (20G44 22D55 57K20)},
  MRNUMBER = {4524205},
MRREVIEWER = {Olga\ Varghese},
       DOI = {10.1112/jlms.12668},
       URL = {https://doi.org/10.1112/jlms.12668},
}

@article {Garland,
    AUTHOR = {Garland, Howard},
     TITLE = {{$p$}-adic curvature and the cohomology of discrete subgroups
              of {$p$}-adic groups},
   JOURNAL = {Ann. of Math. (2)},
  FJOURNAL = {Annals of Mathematics. Second Series},
    VOLUME = {97},
      YEAR = {1973},
     PAGES = {375--423},
      ISSN = {0003-486X},
   MRCLASS = {20J05},
  MRNUMBER = {320180},
MRREVIEWER = {M.\ S.\ Raghunathan},
       DOI = {10.2307/1970829},
       URL = {https://doi.org/10.2307/1970829},
}

@article {EJ,
    AUTHOR = {Ershov, Mikhail and Jaikin-Zapirain, Andrei},
     TITLE = {Property ({T}) for noncommutative universal lattices},
   JOURNAL = {Invent. Math.},
  FJOURNAL = {Inventiones Mathematicae},
    VOLUME = {179},
      YEAR = {2010},
    NUMBER = {2},
     PAGES = {303--347},
      ISSN = {0020-9910,1432-1297},
   MRCLASS = {22D10 (20E08 20E18 20F69)},
  MRNUMBER = {2570119},
MRREVIEWER = {Markus\ Neuhauser},
       DOI = {10.1007/s00222-009-0218-2},
       URL = {https://doi.org/10.1007/s00222-009-0218-2},
}

@book {Marquis_book,
    AUTHOR = {Marquis, Timoth\'ee},
     TITLE = {An introduction to {K}ac-{M}oody groups over fields},
    SERIES = {EMS Textbooks in Mathematics},
 PUBLISHER = {European Mathematical Society (EMS), Z\"urich},
      YEAR = {2018},
     PAGES = {xi+331},
      ISBN = {978-3-03719-187-3},
   MRCLASS = {22-02 (17B99 20E42 20F55 20G44)},
  MRNUMBER = {3838421},
MRREVIEWER = {Ralf\ Koehl},
       DOI = {10.4171/187},
       URL = {https://doi.org/10.4171/187},
}

@unpublished{CapTho,
title={On the {H}owe--{M}oore property for automorphism groups of buildings}, 
      author={Pierre-Emmanuel Caprace and Andreas Thom},
      year={2026},
      note={Preprint, arXiv:2606.27993},
      url={https://arxiv.org/abs/2606.27993}, 
}

@article {GriOpp,
    AUTHOR = {Grinbaum-Reizis, Zohar and Oppenheim, Izhar},
     TITLE = {Curvature criterion for vanishing of group cohomology},
   JOURNAL = {Groups Geom. Dyn.},
  FJOURNAL = {Groups, Geometry, and Dynamics},
    VOLUME = {16},
      YEAR = {2022},
    NUMBER = {3},
     PAGES = {843--862},
      ISSN = {1661-7207,1661-7215},
   MRCLASS = {20J06 (20E42 46C07)},
  MRNUMBER = {4506538},
MRREVIEWER = {Martino\ Garonzi},
       DOI = {10.4171/ggd/676},
       URL = {https://doi.org/10.4171/ggd/676},
}

@article {NoskovVinberg,
    AUTHOR = {Noskov, Guennadi A. and Vinberg, \`Ernest B.},
     TITLE = {Strong {T}its alternative for subgroups of {C}oxeter groups},
   JOURNAL = {J. Lie Theory},
  FJOURNAL = {Journal of Lie Theory},
    VOLUME = {12},
      YEAR = {2002},
    NUMBER = {1},
     PAGES = {259--264},
      ISSN = {0949-5932},
   MRCLASS = {20F55 (20E08 53C23)},
  MRNUMBER = {1885045},
MRREVIEWER = {Daan\ Krammer},
}

@article {CapRAB,
    AUTHOR = {Caprace, Pierre-Emmanuel},
     TITLE = {Automorphism groups of right-angled buildings: simplicity and
              local splittings},
   JOURNAL = {Fund. Math.},
  FJOURNAL = {Fundamenta Mathematicae},
    VOLUME = {224},
      YEAR = {2014},
    NUMBER = {1},
     PAGES = {17--51},
      ISSN = {0016-2736,1730-6329},
   MRCLASS = {20E42 (20E08 20E32 20F65 51E24)},
  MRNUMBER = {3164745},
MRREVIEWER = {Tom\ De Medts},
       DOI = {10.4064/fm224-1-2},
       URL = {https://doi.org/10.4064/fm224-1-2},
}

@book {BH99,
    AUTHOR = {Bridson, Martin R. and Haefliger, Andr\'e},
     TITLE = {Metric spaces of non-positive curvature},
    SERIES = {Grundlehren der mathematischen Wissenschaften [Fundamental
              Principles of Mathematical Sciences]},
    VOLUME = {319},
 PUBLISHER = {Springer-Verlag, Berlin},
      YEAR = {1999},
     PAGES = {xxii+643},
      ISBN = {3-540-64324-9},
   MRCLASS = {53C23 (20F65 53C70 57M07)},
  MRNUMBER = {1744486},
MRREVIEWER = {Athanase\ Papadopoulos},
       DOI = {10.1007/978-3-662-12494-9},
       URL = {https://doi.org/10.1007/978-3-662-12494-9},
}

@article {Rieffel,
    AUTHOR = {Rieffel, Marc A.},
     TITLE = {Induced representations of {$C\sp{\ast} $}-algebras},
   JOURNAL = {Advances in Math.},
  FJOURNAL = {Advances in Mathematics},
    VOLUME = {13},
      YEAR = {1974},
     PAGES = {176--257},
      ISSN = {0001-8708},
   MRCLASS = {46L05 (22D30)},
  MRNUMBER = {353003},
MRREVIEWER = {J.\ W.\ Bunce},
       DOI = {10.1016/0001-8708(74)90068-1},
       URL = {https://doi.org/10.1016/0001-8708(74)90068-1},
}

@article {LarsenPalma,
    AUTHOR = {Larsen, Nadia S. and Palma, Rui},
     TITLE = {Positive definite {$^\ast$}-spherical functions, property
              ({T}) and {$C^\ast$}-completions of {G}elfand pairs},
   JOURNAL = {Math. Proc. Cambridge Philos. Soc.},
  FJOURNAL = {Mathematical Proceedings of the Cambridge Philosophical
              Society},
    VOLUME = {160},
      YEAR = {2016},
    NUMBER = {1},
     PAGES = {77--93},
      ISSN = {0305-0041,1469-8064},
   MRCLASS = {22D25 (22D10)},
  MRNUMBER = {3432330},
MRREVIEWER = {Robert\ S.\ Doran},
       DOI = {10.1017/S0305004115000559},
       URL = {https://doi.org/10.1017/S0305004115000559},
}

@article {BM89,
    AUTHOR = {Barbasch, Dan and Moy, Allen},
     TITLE = {A unitarity criterion for {$p$}-adic groups},
   JOURNAL = {Invent. Math.},
  FJOURNAL = {Inventiones Mathematicae},
    VOLUME = {98},
      YEAR = {1989},
    NUMBER = {1},
     PAGES = {19--37},
      ISSN = {0020-9910,1432-1297},
   MRCLASS = {22E50},
  MRNUMBER = {1010153},
MRREVIEWER = {Marko\ Tadi\'c},
       DOI = {10.1007/BF01388842},
       URL = {https://doi.org/10.1007/BF01388842},
}

@article {BM93,
    AUTHOR = {Barbasch, Dan and Moy, Allen},
     TITLE = {Reduction to real infinitesimal character in affine {H}ecke
              algebras},
   JOURNAL = {J. Amer. Math. Soc.},
  FJOURNAL = {Journal of the American Mathematical Society},
    VOLUME = {6},
      YEAR = {1993},
    NUMBER = {3},
     PAGES = {611--635},
      ISSN = {0894-0347,1088-6834},
   MRCLASS = {22E50},
  MRNUMBER = {1186959},
MRREVIEWER = {David\ Manderscheid},
       DOI = {10.2307/2152779},
       URL = {https://doi.org/10.2307/2152779},
}

@article {Borel,
    AUTHOR = {Borel, Armand},
     TITLE = {Admissible representations of a semi-simple group over a local
              field with vectors fixed under an {I}wahori subgroup},
   JOURNAL = {Invent. Math.},
  FJOURNAL = {Inventiones Mathematicae},
    VOLUME = {35},
      YEAR = {1976},
     PAGES = {233--259},
      ISSN = {0020-9910,1432-1297},
   MRCLASS = {22E50 (20G05)},
  MRNUMBER = {444849},
MRREVIEWER = {Allan\ J.\ Silberger},
       DOI = {10.1007/BF01390139},
       URL = {https://doi.org/10.1007/BF01390139},
}

@book {Matsumoto,
    AUTHOR = {Matsumoto, Hideya},
     TITLE = {Analyse harmonique dans les syst\`emes de {T}its bornologiques
              de type affine},
    SERIES = {Lecture Notes in Mathematics},
    VOLUME = {Vol. 590},
 PUBLISHER = {Springer-Verlag, Berlin-New York},
      YEAR = {1977},
     PAGES = {i+219},
      ISBN = {3-540-08249-2},
   MRCLASS = {22E50 (43A80)},
  MRNUMBER = {579177},
MRREVIEWER = {Yasuhiro\ Asoo},
}

@article {MatsumotoCRAS,
    AUTHOR = {Matsumoto, Hideya},
     TITLE = {Fonctions sph\'eriques sur un groupe semi-simple {$p$}-adique},
   JOURNAL = {C. R. Acad. Sci. Paris S\'er. A-B},
  FJOURNAL = {Comptes Rendus Hebdomadaires des S\'eances de l'Acad\'emie des
              Sciences. S\'eries A et B},
    VOLUME = {269},
      YEAR = {1969},
     PAGES = {A829--A832},
      ISSN = {0151-0509},
   MRCLASS = {22.20 (22.65)},
  MRNUMBER = {263977},
MRREVIEWER = {E.\ Thoma},
}

@article {Ciubotaru,
    AUTHOR = {Ciubotaru, Dan},
     TITLE = {Types and unitary representations of reductive {$p$}-adic
              groups},
   JOURNAL = {Invent. Math.},
  FJOURNAL = {Inventiones Mathematicae},
    VOLUME = {213},
      YEAR = {2018},
    NUMBER = {1},
     PAGES = {237--269},
      ISSN = {0020-9910,1432-1297},
   MRCLASS = {22E50},
  MRNUMBER = {3815566},
MRREVIEWER = {Jeffrey\ D.\ Adler},
       DOI = {10.1007/s00222-018-0790-4},
       URL = {https://doi.org/10.1007/s00222-018-0790-4},
}

@book {Lang,
    AUTHOR = {Lang, Serge},
     TITLE = {{${\rm SL}_2({\bf R})$}},
    SERIES = {Graduate Texts in Mathematics},
    VOLUME = {105},
      NOTE = {Reprint of the 1975 edition},
 PUBLISHER = {Springer-Verlag, New York},
      YEAR = {1985},
     PAGES = {xiv+428},
      ISBN = {0-387-96198-4},
   MRCLASS = {22E45 (11F72)},
  MRNUMBER = {803508},
}

@article {FeitHigman,
    AUTHOR = {Feit, Walter and Higman, Graham},
     TITLE = {The nonexistence of certain generalized polygons},
   JOURNAL = {J. Algebra},
  FJOURNAL = {Journal of Algebra},
    VOLUME = {1},
      YEAR = {1964},
     PAGES = {114--131},
      ISSN = {0021-8693},
   MRCLASS = {20.00},
  MRNUMBER = {170955},
MRREVIEWER = {W.\ E.\ Jenner},
       DOI = {10.1016/0021-8693(64)90028-6},
       URL = {https://doi.org/10.1016/0021-8693(64)90028-6},
}

@article {Duchesne,
    AUTHOR = {Duchesne, Bruno},
     TITLE = {Infinite-dimensional nonpositively curved symmetric spaces of
              finite rank},
   JOURNAL = {Int. Math. Res. Not. IMRN},
  FJOURNAL = {International Mathematics Research Notices. IMRN},
      YEAR = {2013},
            volume={2013},
    NUMBER = {7},
     PAGES = {1578--1627},
      ISSN = {1073-7928,1687-0247},
   MRCLASS = {53C24},
  MRNUMBER = {3044451},
MRREVIEWER = {Zden\v ek\ Du\v sek},
       DOI = {10.1093/imrn/rns093},
       URL = {https://doi.org/10.1093/imrn/rns093},
}

@article {FarHar,
    AUTHOR = {Faraut, Jacques and Harzallah, Kh\'elifa},
     TITLE = {Distances hilbertiennes invariantes sur un espace homog\`ene},
   JOURNAL = {Ann. Inst. Fourier (Grenoble)},
  FJOURNAL = {Universit\'e{} de Grenoble. Annales de l'Institut Fourier},
    VOLUME = {24},
      YEAR = {1974},
    NUMBER = {3},
     PAGES = {xiv, 171--217},
      ISSN = {0373-0956,1777-5310},
   MRCLASS = {43A85},
  MRNUMBER = {365042},
MRREVIEWER = {Wilfried\ Hazod},
       DOI = {10.5802/aif.525},
       URL = {https://doi.org/10.5802/aif.525},
}

@article {BIM,
    AUTHOR = {Burger, Marc and Iozzi, Alessandra and Monod, Nicolas},
     TITLE = {Equivariant embeddings of trees into hyperbolic spaces},
   JOURNAL = {Int. Math. Res. Not.},
  FJOURNAL = {International Mathematics Research Notices},
      YEAR = {2005},
      volume={2005},
    NUMBER = {22},
     PAGES = {1331--1369},
      ISSN = {1073-7928,1687-0247},
   MRCLASS = {20E08 (20F65 57M07)},
  MRNUMBER = {2152540},
       DOI = {10.1155/IMRN.2005.1331},
       URL = {https://doi.org/10.1155/IMRN.2005.1331},
}

@article {MoPy,
    AUTHOR = {Monod, Nicolas and Py, Pierre},
     TITLE = {An exotic deformation of the hyperbolic space},
   JOURNAL = {Amer. J. Math.},
  FJOURNAL = {American Journal of Mathematics},
    VOLUME = {136},
      YEAR = {2014},
    NUMBER = {5},
     PAGES = {1249--1299},
      ISSN = {0002-9327,1080-6377},
   MRCLASS = {53C23 (53C21 57S25)},
  MRNUMBER = {3263898},
MRREVIEWER = {Fernando\ Galaz-Garc\'ia},
       DOI = {10.1353/ajm.2014.0036},
       URL = {https://doi.org/10.1353/ajm.2014.0036},
}

@article {MoPy2,
    AUTHOR = {Monod, Nicolas and Py, Pierre},
     TITLE = {Self-representations of the {M}\"obius group},
   JOURNAL = {Ann. H. Lebesgue},
  FJOURNAL = {Annales Henri Lebesgue},
    VOLUME = {2},
      YEAR = {2019},
     PAGES = {259--280},
      ISSN = {2644-9463},
   MRCLASS = {57S25 (53A35)},
  MRNUMBER = {3978389},
MRREVIEWER = {Leonardo\ Biliotti},
       DOI = {10.5802/ahl.14},
       URL = {https://doi.org/10.5802/ahl.14},
}

@phdthesis{DuchesnePhD, 
  TITLE = {{Des espaces de Hadamard sym{\'e}triques de dimension infinie et de rang fini}},
  AUTHOR = {Duchesne, Bruno},
  URL = {https://theses.hal.science/tel-00673224},
  SCHOOL = {{University of Geneva}},
  YEAR = {2011},
  TYPE = {Theses},
  HAL_ID = {tel-00673224},
  HAL_VERSION = {v1},
}

@book {Robert_book,
    AUTHOR = {Robert, Alain},
     TITLE = {Introduction to the representation theory of compact and
              locally compact groups},
    SERIES = {London Mathematical Society Lecture Note Series},
    VOLUME = {80},
 PUBLISHER = {Cambridge University Press, Cambridge-New York},
      YEAR = {1983},
     PAGES = {ix+205},
      ISBN = {0-521-28975-0},
   MRCLASS = {22D10 (22-01 22D12)},
  MRNUMBER = {690955},
MRREVIEWER = {A.\ Kleppner},
}

@incollection {Gromov,
    AUTHOR = {Gromov, M.},
     TITLE = {Asymptotic invariants of infinite groups},
 BOOKTITLE = {Geometric group theory, {V}ol.\ 2 ({S}ussex, 1991)},
    SERIES = {London Math. Soc. Lecture Note Ser.},
    VOLUME = {182},
     PAGES = {1--295},
 PUBLISHER = {Cambridge Univ. Press, Cambridge},
      YEAR = {1993},
      ISBN = {0-521-44680-5},
   MRCLASS = {20F32 (57M07)},
  MRNUMBER = {1253544},
}

@article {ErshovRall,
    AUTHOR = {Ershov, Mikhail and Rall, Ashley},
     TITLE = {Property ({T}) for {K}ac-{M}oody groups over rings},
   JOURNAL = {J. Algebra},
  FJOURNAL = {Journal of Algebra},
    VOLUME = {500},
      YEAR = {2018},
     PAGES = {221--241},
      ISSN = {0021-8693,1090-266X},
   MRCLASS = {20G44 (20E42 22E40)},
  MRNUMBER = {3765454},
MRREVIEWER = {Walter\ D.\ Freyn},
       DOI = {10.1016/j.jalgebra.2017.01.004},
       URL = {https://doi.org/10.1016/j.jalgebra.2017.01.004},
}

@book {Conway,
    AUTHOR = {Conway, John B.},
     TITLE = {A course in functional analysis},
    SERIES = {Graduate Texts in Mathematics},
    VOLUME = {96},
   EDITION = {Second},
 PUBLISHER = {Springer-Verlag, New York},
      YEAR = {1990},
     PAGES = {xvi+399},
      ISBN = {0-387-97245-5},
   MRCLASS = {46-01 (47-01)},
  MRNUMBER = {1070713},
}

@article {Viola,
    AUTHOR = {Viola, Federico},
     TITLE = {A fixed-point theorem for the action of linear higher-rank
              algebraic groups over local fields on symmetric spaces of
              infinite dimension and finite rank},
   JOURNAL = {Int. Math. Res. Not. IMRN},
  FJOURNAL = {International Mathematics Research Notices. IMRN},
      YEAR = {2026},
            volume = {2026},
    NUMBER = {20},
     PAGES = {Paper No. rnag229},
      ISSN = {1073-7928,1687-0247},
   MRCLASS = {99-06},
  MRNUMBER = {5141849},
       DOI = {10.1093/imrn/rnag229},
       URL = {https://doi.org/10.1093/imrn/rnag229},
}

@article {Duchesne_2023,
    AUTHOR = {Duchesne, Bruno},
     TITLE = {Infinite dimensional representations of orthogonal groups of
              quadratic forms with finite index},
   JOURNAL = {Ann. Fac. Sci. Toulouse Math. (6)},
  FJOURNAL = {Annales de la Facult\'e{} des Sciences de Toulouse.
              Math\'ematiques. S\'erie 6},
    VOLUME = {32},
      YEAR = {2023},
    NUMBER = {2},
     PAGES = {371--396},
      ISSN = {0240-2963,2258-7519},
   MRCLASS = {22E46 (22E65)},
  MRNUMBER = {4631785},
MRREVIEWER = {Mart\'in\ Miglioli},
       DOI = {10.5802/afst.1740},
       URL = {https://doi.org/10.5802/afst.1740},
}

@book {Dixmier,
    AUTHOR = {Dixmier, Jacques},
     TITLE = {Les {$C^*$}-alg\`ebres et leurs repr\'esentations},
    SERIES = {Les Grands Classiques Gauthier-Villars. [Gauthier-Villars
              Great Classics]},
      NOTE = {Reprint of the second (1969) edition},
 PUBLISHER = {\'Editions Jacques Gabay, Paris},
      YEAR = {1996},
     PAGES = {403},
      ISBN = {2-87647-013-6},
   MRCLASS = {46Lxx (01A75 46L10)},
  MRNUMBER = {1452364},
}

@article {CKL,
    AUTHOR = {Caspers, Martijn and Klisse, Mario and Larsen, Nadia S.},
     TITLE = {Graph product {K}hintchine inequalities and {H}ecke {$\rm
              C^\ast$}-algebras: {H}aagerup inequalities, (non)simplicity,
              nuclearity and exactness},
   JOURNAL = {J. Funct. Anal.},
  FJOURNAL = {Journal of Functional Analysis},
    VOLUME = {280},
      YEAR = {2021},
    NUMBER = {1},
     PAGES = {Paper No. 108795, 41},
      ISSN = {0022-1236,1096-0783},
   MRCLASS = {46L07 (46L09 47A63)},
  MRNUMBER = {4156133},
MRREVIEWER = {Adam\ G.\ Skalski},
       DOI = {10.1016/j.jfa.2020.108795},
       URL = {https://doi.org/10.1016/j.jfa.2020.108795},
}

@article {Parkinson,
    AUTHOR = {Parkinson, James},
     TITLE = {Buildings and {H}ecke algebras},
   JOURNAL = {J. Algebra},
  FJOURNAL = {Journal of Algebra},
    VOLUME = {297},
      YEAR = {2006},
    NUMBER = {1},
     PAGES = {1--49},
      ISSN = {0021-8693,1090-266X},
   MRCLASS = {20C08 (05E30 20E42 33C52)},
  MRNUMBER = {2206366},
MRREVIEWER = {Paul-Hermann\ Zieschang},
       DOI = {10.1016/j.jalgebra.2005.08.036},
       URL = {https://doi.org/10.1016/j.jalgebra.2005.08.036},
}

@article {DMSS,
    AUTHOR = {De Medts, Tom and Silva, Ana C. and Struyve, Koen},
     TITLE = {Universal groups for right-angled buildings},
   JOURNAL = {Groups Geom. Dyn.},
  FJOURNAL = {Groups, Geometry, and Dynamics},
    VOLUME = {12},
      YEAR = {2018},
    NUMBER = {1},
     PAGES = {231--287},
      ISSN = {1661-7207,1661-7215},
   MRCLASS = {51E24 (20E08 20E18 20E22 20E32)},
  MRNUMBER = {3781422},
MRREVIEWER = {Guy\ Rousseau},
       DOI = {10.4171/GGD/443},
       URL = {https://doi.org/10.4171/GGD/443},
}

@article {Pont44,
    AUTHOR = {Pontryagin, L. S.},
    TITLE={Hermitian operators in spaces with indefinite metric},
JOURNAL = {Izvestiya Rossiiskoi Akademii Nauk. Seriya Matematicheskaya}, 
    VOLUME = {8},
          YEAR = {1944},
PAGES = {243--280},
}

@unpublished{AnantharamanPopa,
  author = {Anantharaman, Claire and Popa, Sorin},
  title  = {An introduction to {$\mathrm{II}_1$} factors},
  note   = {Book in preparation},
  url    = {https://www.math.ucla.edu/~popa/Books/IIun.pdf}
}

@article {PySanchez,
    AUTHOR = {Py, Pierre and S\'anchez, Arturo},
     TITLE = {Hyperbolic spaces, principal series, and {$\rm O(2,\infty)$}},
   JOURNAL = {Arch. Math. (Basel)},
  FJOURNAL = {Archiv der Mathematik},
    VOLUME = {114},
      YEAR = {2020},
    NUMBER = {1},
     PAGES = {97--105},
      ISSN = {0003-889X,1420-8938},
   MRCLASS = {22D12 (22D20 22E43 22E46 53C35 57T10)},
  MRNUMBER = {4047661},
MRREVIEWER = {Corina\ Ciobotaru},
       DOI = {10.1007/s00013-019-01399-2},
       URL = {https://doi.org/10.1007/s00013-019-01399-2},
}

\end{document}